# An Integrated Optimization Framework for Smart Charging of Electric Bus Fleets under Dynamic Electricity Prices with On-Site Solar Generation, Energy Storage, and V2G operations

Louise Caustur[1], Penelope Hertoghe[1], Tai-Yu Ma[2], Martina Vandebroek[1]

1. KU Leuven University, Naamsestraat 69, 3000 Leuven, Belgium
2. Luxembourg Institute of Socio-Economic Research (LISER), 11 Porte des Sciences, 4366 Esch-sur-Alzette, Luxembourg

## Abstract

The rapid electrification of city bus fleets presents public transportation operators (PTOs) with the complex challenge of managing charging operations to minimize energy costs. Most existing studies on electric bus (EB) charging management rely on a discrete-time-based discretization approach, which is operationally unrealistic and limits their scalability for realistic applications. This study proposes a discrete-event optimization (DEO) approach for daily EB fleet charging management that considers peak power charges, photovoltaic (PV) generation with an energy storage system (ESS), vehicle-to-grid (V2G) operations, and battery degradation costs. We apply the DEO approach to a real-world case in Brussels involving 28 articulated EBs and 232 trips. A set of parametric instances is used to assess computational scalability. The results demonstrate that the DEO formulation can solve instances of realistic size within practical computation times with tight optimality gaps. A thorough cost analysis was conducted to evaluate the added value of V2G benefits and on-site PV generation. Key findings indicate that incorporating demand charges into the optimization reduces daily costs by 5% and decreases the share of peak power costs by 9%, underscoring the importance of load management. Integrating PV and ESS results in a total net cost reduction of up to 56%, with ESS primarily used for energy arbitrage rather than direct bus charging. V2G participation is highly sensitive to battery degradation costs and policy incentives. Combining all extensions results in a 58% reduction in total operational expenses compared to the baseline, demonstrating the significant value of smart (dis)charging tools for PTOs.

## 1. Introduction

The rapid electrification of public bus fleets poses new operational challenges for **public transport operators (PTOs)**, particularly in coordinating **electric bus (EB)** charging operations with limited grid capacity constraints to meet scheduled trip service requirements. Furthermore, fast chargers can trigger high electricity power demand charges if not adequately managed (He et al., 2020a). To address this challenge, PTOs need to implement strategies to reduce demand charges, such as peak shaving, while ensuring the bus service schedule and reliability, which is a critical concern for PTOs. Hence, effective charging schedules are critical to maintaining service reliability and avoiding grid stress. From an operational perspective, electric bus charging system is increasingly coupled with emerging system components to reduce total energy costs. These include time-varying energy prices, peak power charges, heterogeneous charging infrastructure, on-site **photovoltaic (PV)** generation with **energy storage systems (ESS),** and bidirectional **vehicle-to-grid (V2G)** operations. Integrating these emerging system components is essential to reduce daily EB fleet operational costs. However, it brings new challenges to coordinate EB charging decisions with increasing system complexity.

Although numerous studies on electric bus charging optimization have been proposed over the past decades, several limitations remain that compromise the support for cost-efficient charging management of EB fleet. First, most studies consider only a subset of the aforementioned EB charging system components and neglect their interactions. This limits PTO's ability to assess the benefit of coordinated decisions within an integrated optimization framework. Second, most studies formulate the underlying optimization problems based on a **discrete-time-optimization (DTO)** discretization approach, which discretizes the planning horizon into uniform time steps (typically 1 or 5 minutes) and at each time step the model decides whether each EB is charging or not. The main drawback is that charging operations can switch at uncertain times within a given interval, which is unrealistic from an operational standpoint, and that the implied optimization model is very large and has poor scalability for real-world applications. Consequently, the practical application of integrated optimization models, which are usually based on DTO formulation, remains limited for real-world EB networks.

To address these gaps, this study develops a comprehensive and scalable **discrete-event optimization (DEO)** approach that incorporates all the critical system components to support the efficient charging/discharging management for PTOs, thereby reducing their impact on the grid. A DEO-based **Mixed Integer Linear Programming (MILP)** model is proposed that incorporates dynamic electricity tariffs, peak power costs (i.e., demand charges), local PV integration with ESS, a V2G extension, and battery degradation costs. The DEO approach discretizes decision time only at meaningful events (e.g., bus depart/arrival at the depots/terminals or changes in electricity price instead of every minute for the DTO approach), aligning (dis)charging decisions with actual bus operations and dramatically reducing the number of decision variables while providing high-quality solutions. The DEO approach has been successfully applied to conventional EB charging scheduling, ride-pooling, and autonomous dial-a-ride problems recently (Abdelwahed et al., 2020; Gaul et al., 2022; Stallhofer and Parragh, 2025). Using a DEO-based formulation for EB charging scheduling provides a more scalable and realistic representation of charging/discharging operations that account for the complex interactions among different system components and dynamic conditions. The main contributions of this study are listed as follows:

- We develop a DEO-based MILP formulation for the daily charging and discharging scheduling of electric bus fleets. By discretizing time only at relevant events, e.g., bus arrivals or departures at the depots/terminals, electricity price changes, and solar radiation changes, the proposed modeling approach provides operationally realistic solutions and significantly reduces problem size compared to conventional DTO-based approaches. Our results demonstrate that the DEO-based approach can efficiently solve the problem in real-world applications.
- We develop a comprehensive and integrated optimization framework for optimizing the charging schedules of EB networks. The framework jointly captures bidirectional EB charging and discharging, on-site PV generation, and ESS, while accounting for peak demand charges, battery degradation costs, and dynamic electricity prices. This integrated formulation enables PTOs to consistently coordinate complex charging, discharging, and grid interaction decisions in response to dynamic energy prices and on-site PV generation.
- The DEO-based formulation is evaluated on a set of parametric instances with varying fleet sizes and infrastructure characteristics and on a real-world case study for a bus network in Brussels involving 28 EBs and 232 trips with practical considerations. The results show that, within reasonable computational time, near-optimal solutions can be obtained, confirming the scalability and suitability of the DEO-based approach for supporting operational decisions in real-world EB systems.
- A thorough cost analysis is conducted to assess the added value of V2G and on-site renewable generation and to compare it with the baseline scenario. Our results demonstrate that integrating on-site PV generation, ESS, and V2G can reduce EB charging operating costs significantly, with a total decrease of 58% compared to the benchmark. In particular, ESS is used primarily for energy arbitrage rather than for direct EV charging. V2G participation only occurs when the costs of battery degradation and tariff incentives exceed economic thresholds. Through scenario and sensitivity analyses, we explore the impact of battery degradation costs, tariff margins (i.e., margins for injecting energy back into the grid), and green certificates. The analysis provides managerial insights into the trade-off between EV battery replacement costs and V2G revenue streams and between short-term revenue gains and long-term battery wear under different tariff margin scenarios.

The structure of this paper is organized as follows. Section 2 presents a review of relevant literature on electric bus charging optimization, highlighting current modeling approaches and identifying key research gaps. Section 3 presents the modelling framework and formulates the problem as an MILP using the DEO approach. Section 4 describes the realistic case study and its parameter setting based on regular bus lines in Brussels. We consider different modeling scenarios (basic vs. comprehensive) of the EB charging scheduling and compare their performance for the case study. The model scalability is demonstrated via a set of test instances with different fleet sizes and parameter settings. Section 6 conducts a sensitivity analysis to evaluate the impact of varying degradation costs and electricity tariff margins on the system's performance. Finally, conclusions are drawn, and potential directions for future research are discussed.

## 2. Literature review

Research on optimal management and scheduling of electric bus fleets has gained importance in recent years, with a strong focus on infrastructure planning, fleet operation, and charging coordination. This section provides an overview of existing work on the operation of electric bus networks, with particular emphasis on recent extensions (e.g., V2G capability, battery degradation costs, dynamic electricity tariffs, renewable integration) to EB charging scheduling optimization. For a comprehensive review on EB charging scheduling, the reader is referred to Behnia et al. (2024). This section is organized into two distinguished parts: deterministic optimization studies and optimization approaches with uncertainty, followed by a synthesis of research gaps.

### 2.1. Deterministic Models

Early studies on EB charging operation management primarily focused on coordinated charging scheduling under fixed electricity prices and charging infrastructure constraints. More recent works integrate additional options, including dynamic energy prices, V2G operations, battery degradation costs, peak demand charges, and heterogeneous charging technologies. For example, Zheng et al. (2023) developed a MILP model that uniquely addresses daytime and overnight charging scheduling for EBs by integrating grid constraints with dynamic electricity prices. A case study based on a Shanghai bus network shows that the proposed approach achieves a 7.8% cost reduction compared with a first-come, first-served strategy.

Previous studies show that high power demands, especially from fast chargers, can stress the grid and increase peak load charges (Deb et al., 2018; He et al., 2020a). These studies show that peak power appears to be a critical system cost, particularly when fast chargers are installed for EB charging operations. He et al. (2020b) propose a model to minimize total charging costs by optimizing the scheduling and management of a fast-charging electric bus system across multiple lines. The model incorporates opportunity charging, demand charges, and time-of-use (TOU) electricity rates, and is particularly insightful for understanding the factors influencing charging decisions and their impact on cost optimization.

As noted earlier, a central modeling challenge in EB charging scheduling is that most existing formulations rely on discrete-time discretization schemes, which are often inconsistent with practical operational constraints (Cheng et al., 2019; He et al., 2020a, 2020b; Manzolli et al., 2022a, b, 2024, 2025; Liu et al., 2024; Son et al., 2025). These studies use uniform discrete-time discretization with a time interval of 5 to 60 minutes to ensure tractability. Furthermore, the DTO-based approach can lead to scalability issues when coordinating charging schedules that involve complex interactions among different system components in real-world applications. Abdelwahed et al. (2020) provide a systematic comparison between the DTO-based and DEO-based approaches for EB fleet charging management. In that study, the DTO approach discretizes the planning horizon with fine-grained time steps (1 minute), whereas the DEO discretizes time only with respect to relevant events (e.g., a bus arriving or departing at a terminal/depot, an electricity price change, or the end of the planning horizon). Additional constraints are introduced for the DTO approach to ensure operational consistency across consecutive time steps. The DEO-based MILP shows significant gains in computational performance and scalability, with good solution quality (<0.3% difference in objective values), thanks to its substantial reduction in problem size (number of decision variables). However, their study failed to consider V2G operations and on-site renewable generation and peak power costs.

Another research stream incorporates V2G operations and battery degradation costs into EV charging management. For example, Manzolli et al. (2022b) develop a MILP model minimizing the daily operating costs to help the PTO control the energy required to operate a bus fleet while accounting for V2G interactions with the grid. The results demonstrate that the decision to engage in V2G operations depends on the trade-off between the battery replacement costs and energy selling prices. In another study (Manzolli et al., 2022a), the authors propose a MILP model to coordinate the charging events of an electric bus fleet system under a semi-empirical battery degradation model. They consider battery ageing costs due to cycle ageing (i.e., costs associated with frequent charging and discharging of the batteries). The effect of weather conditions on battery degradation is also considered by looking at a winter scenario and a summer scenario. The findings of both studies do not allow for V2G operations because the battery degradation costs are higher than the electricity selling prices. However, the issue of whether the current decreasing trend for battery replacement costs will incentivize the use of V2G in public transport systems is still underexplored. More recent work by Manzolli et al. (2025) integrates TOU energy prices, V2G operation, and battery degradation costs within a robust optimization framework. While this study provides valuable insights into the potential gains of coordinated charging under uncertainty, their approach remains computational limited (largest computational instance is limited to 10 buses and 10 trips) due to their reliance on DTO-based MILP formulations.

Integrating renewable energy has become a critical focus in the effort toward carbon neutrality for the electrification of EB networks. The installation of local PV panels enables PTOs to charge their buses at zero cost when solar radiation is available. Additionally, it contributes to the reduction of carbon emissions, alleviates street-level air pollution, and improves on-site solar energy consumption (Ren et al., 2022). PV panel installations are often coupled with an ESS, typically battery-based, that can help reduce the unstable supply of solar-generated electricity (Dougier et al., 2023). These systems also enable to store surplus energy generated during periods of high solar availability, which can later be used to recharge the buses. In addition to supporting charging operations, the ESS can also engage in the market and export electricity back to the grid. This dual role makes it a valuable complement, not only to satisfy the charging demand of the EBs but also to reduce charging costs through electricity price arbitrage, also referred to as energy arbitrage (Ding et al., 2015). Cheng et al. (2019) developed an aggregation strategy to coordinate the operations of a bus network consisting of multiple lines integrated with on-site PV generation and a stationary ESS. The authors apply their model to a case study and find that implementing an aggregation strategy can reduce the total cost for the PTO by 17% compared to a setup where each bus line has its own charging station. Furthermore, it shows that with an ESS and optimized PV integration, total costs can be decreased in comparison with no ESS. Similarly, Liu et al. (2024) proposed an EB charging scheduling model by integrating on-site PV generation and ESS. The results show a reduction in daily costs by 15% to 33%, depending on weather conditions, highlighting the financial benefits of renewable energy integration. However, these models failed to integrate V2G operations and battery ageing costs, thereby overlooking the potential cost savings of V2G participations.

### 2.2. Models with Uncertainty

Various sources of uncertainty affect real-world EB operations, including energy consumption, travel times, renewable energy generation, and operational constraints. The challenge of handling uncertainty in the EB network has led to various modeling approaches to tackle this uncertainty, ranging from robust optimization to stochastic programming.

Hu et al. (2022) look at the implementation of fast-charging stations along the bus trajectory to allow for en-route charging to reduce the battery sizes of the buses. The paper aims to locate the fast chargers at selected bus stops and provide an optimal corresponding charging schedule while considering the stochasticity of travel time and of passengers' boarding and alighting time. The results highlight that en-route fast charging is effective in maintaining operational schedules while reducing energy consumption costs compared to conventional combustion buses. Additionally, Zhuang and Liang (2021) propose a stochastic energy management method for electric bus charging stations using a distributionally robust Markov decision process. The model incorporates PV energy generation, integrated ESS, and V2G capabilities to create cost-effective charging strategies. Case studies from St. Albert (Canada) and IEEE test feeders demonstrate significant cost savings and reduced grid impact, highlighting the efficacy of the method in addressing uncertainty in bus loads and PV generation. In addition to robust and stochastic methods, Abdelwahed et al. (2024) propose a real-time decision support system to optimize charging schedules and mitigate the impact of operational uncertainties. Their framework addresses uncertainties such as operational delays, energy consumption variability, and the intermittency of renewable energy resources by leveraging real-time data, predictions, and mathematical optimization. Key insights include the importance of optimizing charging schedules, integrating renewable energy, and recognizing that opportunity fast-charging bus networks are highly sensitive to operational uncertainties. The paper suggests that ESSs could further attenuate the variability associated with renewable energy integration. However, its modeling approach did not account for EBs with V2G capability or peak power costs.

Finally, Manzolli et al. (2025) detail and aggregate several key characteristics discussed in previous works. They propose a robust optimization model that accounts for uncertainty about energy consumption and factors such as battery aging, TOU tariffs, V2G integration, and operational constraints such as fixed schedules. Through a case study in Coimbra, a medium-sized city in Portugal, the findings indicate that a significant cost reduction can be achieved by using coordinated charging: 37% for the deterministic model and 12% for the robust model compared to usual charging scenarios. However, the problem is formulated as a MILP based on the DTO approach, which limits its applications to a small bus network with 10 buses and 10 trips.

The above literature review reveals three research gaps. First, no existing study jointly integrates dynamic energy prices, capacity-based demand charges, on-site PV generation, ESS, V2G, and battery degradation costs within an integrated DEO-based optimization framework for EB (dis)charging management. Second, despite the scalability benefits and enhanced operational realism offered by the DEO-based approach, most existing studies rely on DTO-based formulations, thereby limiting their practical applicability in real-world settings. Third, although incorporating multiple sources of uncertainty can enhance operational robustness and service reliability, day-ahead charging scheduling decision-support tools remain essential to ensure operational feasibility and cost efficiency in practice**.** PTOs can subsequently extend such deterministic day-ahead scheduling models into real-time decision-support systems by dynamically adapting planned charging schedules in response to operational disruptions, demand fluctuations, and energy-system uncertainties.

Table 1 summarizes the range of charging methods and system components of existing studies, as well as the differences and new contributions of this study.

Table 1. Summary of charging methods and optimization techniques.

| **Paper** | **Opportunity Charging** | **Overnight Charging** | **Uncertain Energy consump. of EB** | **Renewable Energy Integration** | **Dynamic electricity price** | **V2G** | **Battery Degradation** | **Peak Power cost** | **DEO/ DTO** |
|---|---|---|---|---|---|---|---|---|---|
| Cheng et al. (2019) | | ✓ | | ✓ | ✓ | | | | DTO |
| Abdelwahed et al. (2020) | ✓ | ✓ | | | ✓ | | | | DEO+ DTO |
| Zhuang & Liang, (2021) | ✓ | ✓ | ✓ | ✓ | ✓ | ✓ | ✓ | | DTO |
| Ren et al. (2022) | ✓ | ✓ | | | | | | | DTO |
| Liu et al. (2022) | | ✓ | ✓ | | | | | | DTO |
| Manzolli et al. (2022a) | ✓ | ✓ | | | ✓ | ✓ | ✓ | | DTO |
| Manzolli et al. (2022b) | ✓ | ✓ | | | ✓ | | | ✓ | DTO |
| Liu et al. (2024) | ✓ | ✓ | | | ✓ | | | ✓ | DTO |
| Manzolli et al. (2024) | ✓ | ✓ | | | ✓ | ✓ | ✓ | ✓ | DTO |
| Abdelwahed et al. (2024) | ✓ | ✓ | ✓ | ✓ | ✓ | | | | DEO |
| Manzolli et al. (2025) | | ✓ | ✓ | | ✓ | ✓ | ✓ | ✓ | DTO |
| Son et al. (2025) | | ✓ | ✓ | | ✓ | ✓ | | | DTO |
| **This paper** | ✓ | ✓ | | ✓ | ✓ | ✓ | ✓ | ✓ | **DEO** |

## 3. Methodology

***Nomenclature***

***Sets***

| | |
|---|---|
| $K$ | Set of electric buses |
| $E$ | Set of events in the network with index *e*, each time slot of event e lies between events *e* and *e + 1* |
| $J$ | Set of depots |
| $N_j$ | Set of chargers located at depot $j \in J$ |
| $I$ | Set of trips |
| $L$ | Set of discrete power levels |

***Indices***

| | |
|---|---|
| $k \in K$ | Index for an electric bus |
| $e \in E$ | Index for event *e* |
| $j \in J$ | Index for a bus depot |
| $n \in N_j$ | Index for a charger located at bus depot $j \in J$ |
| $i \in I$ | Index for a trip |
| $l \in L$ | Index for a power level |

***Parameters***

| | |
|---|---|
| $\gamma_i$ | Average energy consumption for trip *i* (kWh/km) |
| $\eta_{j,n}^{char}, \eta_{j,n}^{dis}$ | Charging/discharging efficiency of charger *n* at depot *j* (%) |
| $\alpha_{j,n}$ | Charging power of charger *n* at depot *j* (kW) |
| $\beta_{j,n}$ | Discharging power of charger *n* at depot *j* (kW) |
| $T_e$ | Time at which event *e* occurs and time slot of event *e* starts (minute) |
| $\Delta T_e$ | Length of time slot of event *e* ($T_{e+1} - T_e$)(minute) |
| $C_k^{bat}$ | Total capacity of the battery of bus *k* (kWh) |
| $E_k^{min}, E_k^{max}$ | Minimum/Maximum SOC allowed for bus *k* (%) |
| $E_k^0$ | Initial SOC of bus *k* (%) |
| $E_k^{end}$ | Minimum SOC maintained for bus *k* at the end of an operations' day (%) |
| $\overline{H}_j$ | Total storage capacity of the ESS at depot *j* (kWh) |
| $SOC^{min}$ | Minimum SOC allowed for the batteries at ESS (%) |
| $U_l^{pow}$ | Power at level *l* (kW) |
| $U^{max}$ | Maximum total contracted power level (kW) |
| $U_l^{price}$ | Purchasing price for power level *l* (€/day) |
| $N_k^{cyc}$ | Maximum number of cycles that the battery of bus *k* can last |
| $R_k$ | Battery replacement cost of battery of bus *k* (€/kWh) |
| $\rho_e^+, \rho_e^-$ | Electricity purchasing/selling price during time slot of event *e* (€/kWh) |
| $\tau^m$ | Minimum charging/discharging time required to increase lifespan of the battery |
| $Q_{j,e}$ | Solar PV electricity yield at depot *j* during time slot of event *e* (kWh) |
| $N_j$ | Number of chargers present at depot *j* |
| $\theta_e$ | Binary parameter indicating if PV panels generate electricity during time slot of event *e* |
| $V_e$ | Minimum number of time slots required to charge/discharge if a bus starts charging/discharging at occurrence of event *e* to meet the minimum charging time requirement |
| $M_{j,e}$ | Binary parameter indicating if there is a bus arrival or departure at depot *j* at the occurrence of event *e* |
| $t_i$ | Average speed during the trip *i* (km/h) |
| $\epsilon_k^g$ | Index of time slot at which bus *k* arrives after finishing all trips |
| $l_{j,k,e}$ | Binary parameter indicating if bus *k* is present at depot *j* during time slot of event *e* |
| $b_{k,i,e}$ | Binary parameter indicating if bus *k* serves trip *i* during time slot of event *e* |
| $f_k$ | Total energy taken from battery of bus *k* through its lifespan (kWh) |
| $\epsilon^f$ | Last event of the day before the next operations' day starts (last minute of the 24hr period) |
| ***Decision Variables*** | |
| $x_{j,k,n,e}$ | Binary variable indicating if bus *k* is using charger *n* at depot *j* during time slot of event *e* to charge |
| $\delta_{j,k,n,e}^{char}$,$\delta_{j,k,n,e}^{dis}$ | Charging/discharging duration of bus *k* using charger *n* at depot *j* during time slot of event *e* |
| $\lambda_{k,e}^{char}$,$\lambda_{k,e}^{dis}$ | Binary variable indicating if bus *k* started charging/discharging at occurrence of event *e* |
| $v_{k,e}^{char}$,$v_{k,e}^{dis}$ | Binary variable that equals 1 when $\sum_{j\in J}\sum_{n\in N_j} x_{j,k,n,e}\, x_{j,k,n,e-1} = 1$ , if bus *k* is charging/discharging during both time slot of event *e*−1 and time slot of event *e* (i.e., if a continuous charging session is active at event *e*) |
| $y_{j,k,n,e}$ | Binary variable indicating if bus *k* is using charger *n* at depot *j* during time slot of event *e* to discharge |
| $u_l^{peak}$ | Binary variable indicating if power level *l* is taken as peak power |
| $e_{k,e}$ | Energy level of bus k at the occurrence of event e (kWh) |
| $w_{j,e}^{+}$,$w_{j,e}^{-}$ | Electricity purchased/sold by depot *j* during time slot of event *e* (kWh) |
| $d_{k,e}$ | Total degradation cost of battery of bus *k* at occurrence of event *e* (€) |
| $\mu_{j,e}$ | Usage of solar energy for charging buses at depot *j* during time slot of event *e* when the PV panels generate electricity (kWh) |
| $z_{j,e}$ | Usage of solar energy for charging buses at depot *j* during time slot of event *e* when the PV panels do not generate electricity (kWh) |
| $\pi_{j,e}$ | Amount of solar energy stored a bus depot *j* during time slot of event *e* when PV panels generate electricity (kWh) |
| $H_{j,e}$ | Remaining electricity at the ESS of bus depot *j* at occurrence of event *e* (kWh) |
| $a_{k,e}^{char}$,$a_{k,e}^{dis}$ | Binary variable indicating if bus *k* stops charging/discharging during time slot of event *e* |
| $m_{j,e}^{anc}$ | Electricity reinjected into the grid from the ESS during time slot of event *e* at depot *j* |

3.1. Problem description

We consider a deterministic V2G bus charging scheduling optimization problem with local PV installation and ESS to meet every scheduled trip of a timetabled bus network. The objective of PTOs is to minimize the daily operating cost. The operating cost comprises four monetary components: (i) energy purchased from the grid, (ii) a capacity-based peak power charge, (iii) EB battery-degradation penalty, and (iv) revenues earned by discharging either bus batteries (V2G) or an ESS from PV-produced energy into the grid. Charging infrastructure (i.e., number of chargers, type, location) is treated as a fixed strategic input: these assets are decided once during network planning and remain unchanged in day-to-day operation. Each charger is modeled as an independent unit with its own rated power and efficiency, so heterogeneous equipment can be represented explicitly. Scheduling of buses is treated as input data which brings this study to optimize and decide for every layover whether a bus should charge, discharge, or remain idle, and at what power level, so that total daily cost is minimized. Input data supplied by the PTO each day includes the day-ahead TOU electricity tariff, PV generation forecasts, and the initial state-of-charge of every battery. By jointly optimizing grid imports, PV self-consumption, EB dispatch and V2G participation, the model balances the PTO's preference for operationally convenient charging windows, often coinciding with high tariff periods, against its need to limit charging costs and peak load. Note that the proposed approach does not explicitly consider battery degradation costs associated with ESS charging and discharging operations due to limited depot-specific cycling/lifetime data. Consequently, the operational benefits of the ESS may be slightly overestimated and can be considered as an upper bound on its economic contribution. We mentioned this limitation in Section 6.

A MILP model is developed to model the EB charging scheduling problem with tailored extensions. By embedding different sources of energy, degradation cost, peak-power pricing, heterogeneous charger technologies, and multiple energy sources in one MILP, the proposed framework unifies elements that have so far been studied only in isolation and offers PTOs a practical decision tool for the future of electric bus fleets. Given the complexity of the problem, we develop a DEO-based MILP formulation for its advantage in scalability. The 24-hr planning horizon is divided into time slots according to relevant events that occur throughout the planning horizon. For this model, **events** are defined as: (i) a bus arrives at a depot/terminal; (ii) a bus leaves a depot/terminal; (iii) electricity price changes; (iv) solar radiation changes; (v) the last minute of the 24-hr period. All these events mark the start of a new time slot. An ***event*** is defined as a point in time, whereas a ***time slot*** is a period of time with a fixed duration. Simultaneous events are merged into a single event index. The start and the end of a time slot are defined by events. The time granularity of the data is in minutes, as bus schedules are usually indicated up to minutes in terms of time units. To stay consistent with the syntax used in the implementation, the first event takes index 1. It is also useful to note that the term depot is used to refer to both terminals and overnight depots. In this context, a depot indicates an infrastructure that contains chargers. Figure 1 shows how events are defined and how time slots and events are indexed.

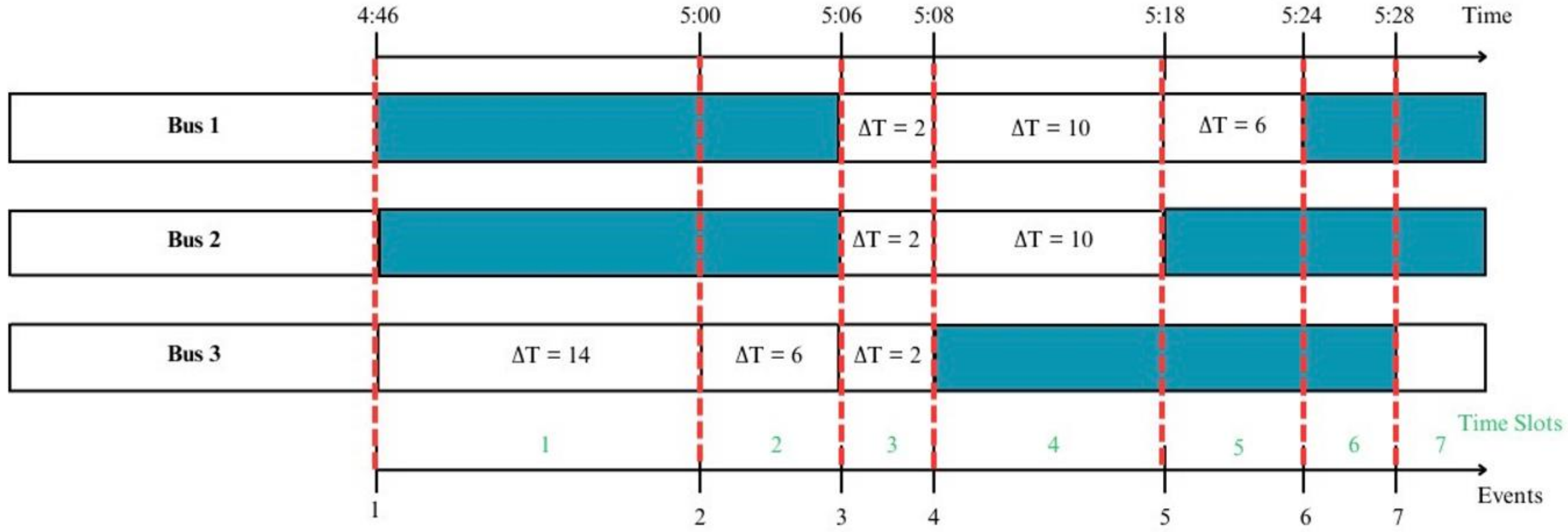

Figure 1. Definition of events and time slots. The shaded areas represent the time slots during which the buses are unavailable for charging or discharging because they are serving a trip. Note that at 5:00, no arrival or departure of buses defines the occurrence of the event, but the event still exists because there is a change in the electricity price.

3.2. Assumptions

The assumptions of the proposed modeling approach are summarized as follows.

a. All operational inputs are known deterministically at the day-ahead planning stage, including bus trip assignment and timetables, TOU energy price, PV generation forecast, initial state-of-charge of EBs. Uncertainty in travel time, energy consumption, PV generation, and electricity prices is not explicitly modelled. These assumptions allow us to assess the structural contribution of each system component.

b. Each charger is modeled as an independent resource with constant power and efficiency, allowing accommodating heterogeneous charging infrastructure configurations
c. Bus energy consumption model is deterministic as a function of average speed. We use a calibrated EB energy consumption function to estimate the energy consumption of each bus trip (see Section 4). Although this approach ignores the influence of factors such as weather, passenger loads, and road slope on EB energy consumption, it is adequate for deriving a system-level EB charging schedule and ensuring the tractability of the problem.
d. Battery charge and discharge are assumed linear with fixed charging and discharging efficiencies. EB battery degradation cost is modeled as a function of the battery capacity, the replacement cost of the battery, and the total number of lifecycles that a battery can last. Degradation costs are considered only when they result from discharging activities, as the model is designed to assess whether participating in V2G operations remains economically attractive despite the added degradation burden.
e. The battery degradation of ESS system is not considered in this study due to the lack of reliable depot-specific cycling and lifetime data. This assumption will slightly overestimate the economic value of ESS component. We acknowledge this limitation in the sensitivity analysis and discussion.
f. Each charging or discharging session must satisfy a minimum duration to reflect practical operational constraints. Furthermore, we assume minimum state-of-charge needs to be maintained for each EB and ESS.

3.3. Mathematical formulation

The integrated EB charging management framework is provided in Figure 2. EBs can charge or discharge when connected to a charger. When the PV panels generate energy during the time slot of event $e$ at depot $j$, it can be either used for EB charging ($\mu_{j,e}$) or stored in ESS ($\pi_{j,e}$). ESS can send energy for charging EBs using solar energy when the PV panels do not generate electricity ($z_{j,e}$) or sell back the energy into the grid ($m_{j,e}^{anc}$) for ancillary service. Note that solar energy generation is based on day-ahead predictions (deterministic). Adapting uncertain renewable generation to support real-time decisions of (dis)charging operations could be a future extension of the present work. The system also determines how much energy to sell to or buy from the grid for each event $e$ and depot $j$. The nomenclature of the model can be found at the beginning of this section, where all sets, indices, parameters, and decision variables can be found with their respective definitions.

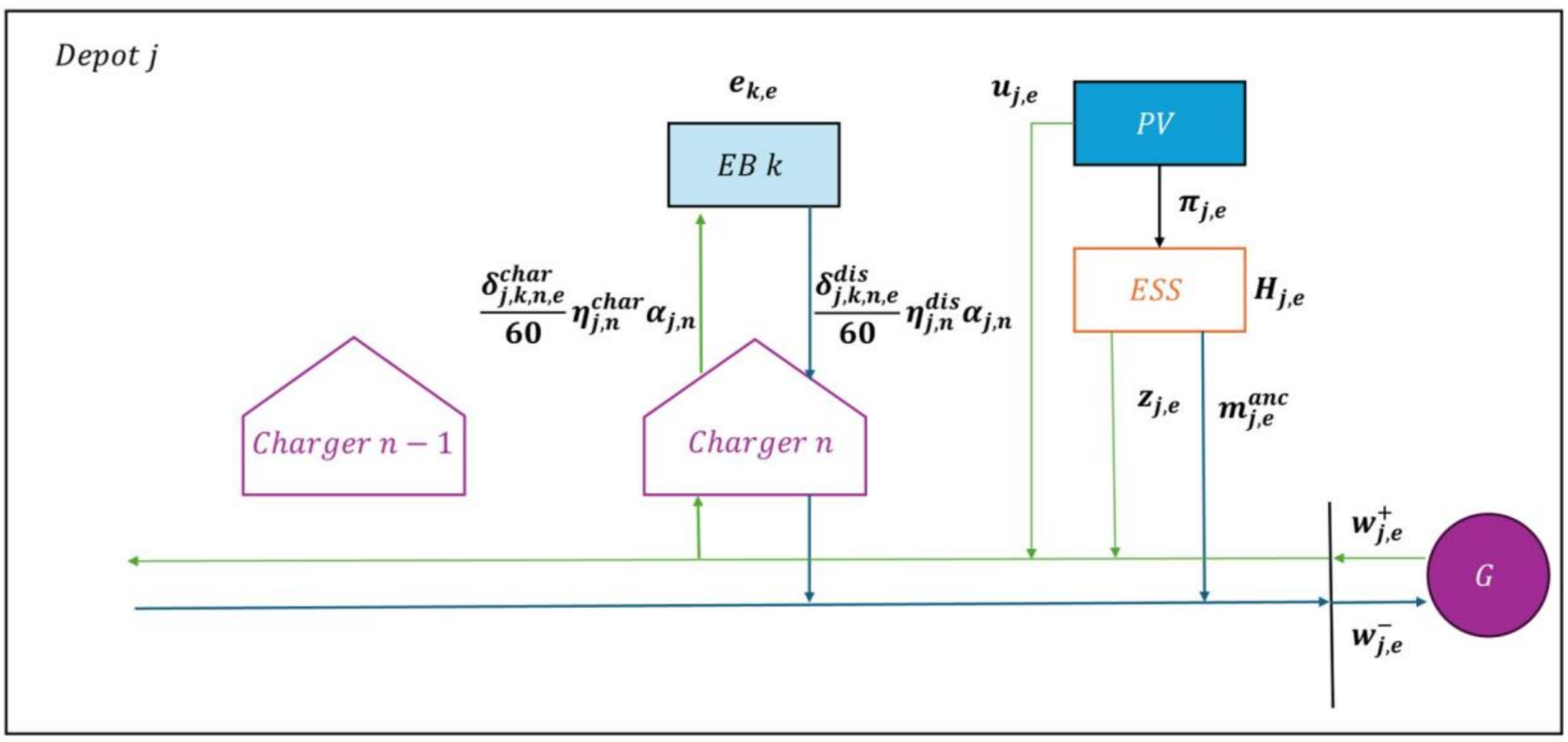

Figure 2. Framework of the DEO model with key decision variables characterizing energy flows between EBs, chargers, PV panels, ESS, and the grid (G).

We detail the DEO-based MILP model (also referred to as the DEO model) as follows.

**Objective function**

$$\min \sum_{l \in L} U_l^{\text{price}} u_l^{\text{peak}} + \sum_{j \in J} \sum_{e \in E} \rho_e^{+} w_{j,e}^{+} - \sum_{j \in J} \sum_{e \in E} \rho_e^{-} w_{j,e}^{-} + \sum_{k \in K} \sum_{e \in E} d_{k,e} \tag{1}$$

The objective function (1) minimizes the overall costs over a 24-hour period if set $E$ (events) covers this time frame. The first cost component is the peak power cost which was incorporated to control the impact on the grid and avoid grid overload. $L$ is the set of peak power levels, while $u_l^{\text{peak}}$ and $U_l^{\text{price}}$ correspond to the power at level $l$ and the purchasing price for power level $l$, respectively. The second and third cost components are the cost of energy purchased from the grid and the revenue generated from selling energy back to the grid, thanks to V2G technology and energy arbitrage. $w_{j,e}^{+}$ ($w_{j,e}^{-}$) denotes the electricity purchased (sold) by depot (or terminal) $j$ during the time slot of event $e$. Finally, battery degradation costs due to discharging activities are taken into account by summing over all the degradation costs incurred by the buses over a day. $d_{k,e}$ denotes the battery degradation costs of bus $k$ at the occurrence of event $e$.

**Charging and Discharging Constraints**

$$\sum_{j \in J} \sum_{n \in N_j} x_{j,k,n,e} + \sum_{j \in J} \sum_{n \in N_j} y_{j,k,n,e} \leq 1 - \sum_{i \in I} b_{k,i,e} \;\; \forall k \in K, e \in E \tag{2}$$

Constraint (2) ensures that a bus $k$ occupies one and only one charger $n$ during time slot of event $e$ to charge or discharge. Additionally, it prevents the buses to either charge ($x_{j,k,n,e} = 1$) or discharge ($y_{j,k,n,e} = 1$) whilst they are serving a trip ($b_{k,i,e} = 1$).

$$\sum_{k \in K} x_{j,k,n,e} + \sum_{k \in K} y_{j,k,n,e} \leq 1 \;\; \forall j \in J, \forall n \in N_j, e \in E \tag{3}$$

Constraint (3) ensures that at most one bus can occupy charger $n$ at depot $j$ during time slot of event $e$.

$$\sum_{k \in K} x_{j,k,n,e} + \sum_{k \in K} y_{j,k,n,e} \leq 1 \;\; \forall j \in J, \forall n \in N_j, e \in E \tag{4}$$

Constraint (4) makes sure that a bus can only engage in a charging or discharging activity at a given depot $j$ if it is present at that depot during time slot of event $e$.

$$e_{k,e+1} = e_{k,e} - \sum_{i \in I} \gamma_i t_i b_{k,i,e} \frac{\Delta T_e}{60} + \sum_{j \in J} \sum_{n \in N_j} \frac{\delta_{j,k,n,e}^{char}}{60} \eta_{j,n}^{char} \alpha_{j,n} - \sum_{j \in J} \sum_{n \in N_j} \frac{\delta_{j,k,n,e}^{dis}}{60} \frac{1}{\eta_{j,n}^{dis}} \beta_{j,n} \quad \forall k \in K, e \in E \setminus \epsilon^f \tag{5}$$

Constraint (5) keeps track of the energy level of bus $k$ at the occurrence of event $e$. The energy level at event $e + 1$ is defined as the energy level at the previous event, minus the energy consumed for a trip during the previous time slot, plus the energy charged, and minus the energy discharged to the grid, all during the same time slot. Note that a bus can engage in at most one of these three activities during any time slot.

$$\sum_{n \in N_j} \sum_{k \in K} \frac{\delta_{j,k,n,e}^{char}}{60} \eta_{j,n}^{char} \alpha_{j,n} = w_{j,e}^{+} + \mu_{j,e} + z_{j,e} \;\; \forall j \in J, e \in E \tag{6}$$

$$\sum_{n \in N_j} \sum_{k \in K} \frac{\delta_{j,k,n,e}^{dis}}{60} \frac{1}{\eta_{j,n}^{dis}} \beta_{j,n} + m_{j,e}^{anc} = w_{j,e}^{-} \ \forall j \in J, e \in E \tag{7}$$

Constraints (6) and (7) compute the amount of energy purchased and sold to the grid during time slot of event $e$. The amount of energy purchased ($w_{j,e}^{+}$) takes into account the energy drawn from PV panels ($\mu_{j,e}$) and the energy drawn from the ESS ($z_{j,e}$). The amount of energy sold to the grid ($w_{j,e}^{-}$) includes the energy reinjected from the ESS ($m_{j,e}^{anc}$) and from the buses.

**Energy level constraints for batteries of buses**

$$e_{k,e} \geq C_k^{bat} E_k^{\min}, \forall k \in K, e \in E \tag{8}$$

$$E_k^{\max} C_k^{bat} \geq e_{k,e}, \forall k \in K, e \in E \tag{9}$$

Constraints (8) prevents the energy in the battery from dropping below a certain threshold whilst constraint (9) prevents it from being above the battery's maximum allowed capacity. These boundaries are set to maximize the life expectancy of the batteries, as charging them or discharging them entirely negatively affects their performance capacities.

$$e_{k,1} = C_k^{bat} E_k^1, \forall k \in K \tag{10}$$

Constraint (10) initializes the battery level at the first event of the operational day (i.e.: $e = 1$ ).

$$e_{k,\epsilon^f} \geq C_k^{bat} E_k^{end}, \forall k \in K \tag{11}$$

Constraint (11) enforces a minimum energy level at the last event of the operational day to ensure that the fleet begins the next day with sufficient energy to operate.

$$H_{j,e} \leq \bar{H}_j \ \forall j \in J, e \in E \tag{12}$$

$$H_{j,e} \geq SOC_{\min} \cdot \bar{H}_j \ \forall j \in J, e \in E \tag{13}$$

Constraints (12) and (13) limit the energy level of the ESS ( $H_{j,e}$ ) with an upper and a lower bound to maximize the ESS's performance over its lifetime.

$$H_{j,e+1} = H_{j,e} + \pi_{j,e} \cdot \theta_e - z_{j,e} \cdot (1 - \theta_e) - m_{j,e}^{anc} \ \forall j \in J, e \in E \setminus \epsilon^f \tag{14}$$

Constraint (14) tracks the energy level of the ESS at occurrence of event $e$. The energy level of the ESS at the occurrence of event $e$ is defined as its energy level at occurrence of the previous event. If PV panels generate electricity (i.e., $\theta_e = 1$) during previous time slot, then solar energy stored ($\pi_{j,e}$) during previous time slot is added. If there is no sun (i.e., $\theta_e = 0$) during the previous time slot, then the energy used to recharge the buses ($z_{j,e}$) during previous time slot is subtracted. Finally, the energy re-injected into the grid for energy arbitrage ($m_{j,e}^{anc}$) during previous time slot is also subtracted.

$$m_{j,e}^{anc} \leq \bar{H}_j, \forall j \in J, e \in E \tag{15}$$

Constraint (15) limits the amount of electricity re-injected into the grid by the ESS during time slot of event $e$ to the maximum capacity of the ESS.

$$H_{j,\epsilon} f \geq SOC_{\min} \bar{H}_j, \forall j \in J \tag{16}$$

Constraint (16) sets the energy level in the ESS during the last event of the operations' day at a minimum level.

$$m_{j,\epsilon^f}^{anc} = 0, \forall j \in J \tag{17}$$

Constraint (17) prohibits the ESS from injecting electricity into the grid during the last time slot of the day because it represents the transition between 2 days.

$$H_{j,1} = SOC_{\min}\bar{H}_j, \forall j \in J \tag{18}$$

Constraint (18) initializes the energy level in the ESS at the occurrence of the first event of the day ($e = 1$).

**Battery degradation Constraints**

$$d_{k,e} = \frac{R_k C_k^{bat}}{f_k} \sum_{j \in J} \sum_{n \in N_j} \beta_{j,n} y_{j,k,n,e} \frac{\Delta T_e}{60}, \forall k \in K, e \in E \tag{19}$$

Constraint (19) defines the total degradation cost of battery of bus $k$ at event $e$ as a function of the battery capacity, the replacement cost of the battery and the total number of lifecycles that a battery can last. $R_k$ is the battery replacement cost of battery of bus k. $C_k^{bat}$ is total capacity of the battery of bus k. $f_k$ is the total energy taken from battery of bus k through its lifespan.

**Charging Duration Determination**

$$\delta_{j,k,n,e}^{char} \geq 1 - \Delta T_e \left(1 - x_{j,k,n,e}\right) \forall k \in K, e \in E, j \in J, n \in N_j \tag{20}$$

$$\delta_{j,k,n,e}^{char} \leq \Delta T_e x_{j,k,n,e} \forall k \in K, e \in E, j \in J, n \in N_j \tag{21}$$

Constraints (20) and (21) make the link between the binary variable $x_{j,k,n,e}$ that indicates whether a bus is charging or not at a given depot and charger during time slot of event $e$ and the charging duration variable $\delta_{j,k,n,e}^{char}$. A big-M constraint is used where the big-M is set to the smallest feasible value (i.e., $M_e = \Delta T_e$ ).

$$\delta_{j,k,n,e}^{dis} \geq 1 - \Delta T_e \left(1 - y_{j,k,n,e}\right), \forall k \in K, e \in E, j \in J, n \in N_j \tag{22}$$

$$\delta_{j,k,n,e}^{dis} \leq \Delta T_e y_{j,k,n,e}, \forall k \in K, e \in E, j \in J, n \in N_j \tag{23}$$

Constraints (22) and (23) make the link between the binary variable $y_{j,k,n,e}$ that indicates whether a bus is discharging or not at a given depot and charger during time slot of event $e$ and the discharging duration variable.

**Number of chargers constraint**

$$\sum_{n \in N_j} \sum_{k \in K} \left(x_{j,k,n,e} + y_{j,k,n,e}\right) l_{j,k,e} \leq \left|N_j\right|, \forall j \in J, e \in E \tag{24}$$

Constraint (24) ensures that the number of buses charging or discharging at a given depot for any event does not exceed the number of chargers $|N_j|$ available at depot $j$.

**Minimum Charging Time and Minimum Discharging Time**

$$\sum_{j\in J}\sum_{n\in N_j}\sum_{e^*=e}^{e+V_e-1}\delta_{j,k,n,e^*}^{char} \geq \tau^m \cdot \lambda_{k,e}^{char}, \forall k \in K, e \in E \quad (25)$$

Constraint (25) makes sure that if bus $k$ starts charging at event $e$ (i.e., $\lambda_{k,e}^{\text{char}} = 1$), then its charging duration is above the minimum charging time $\tau^m$.

$$\sum_{j\in J}\sum_{n\in N_j}\sum_{e^*=e}^{e+V_e-1} x_{j,k,n,e^*} \geq V_e\lambda_{k,e}^{char}, \forall k \in K, e \in E \quad (26)$$

Constraint (26) makes sure that if bus $k$ starts charging at event $e$ (i.e., $\lambda_{k,e}^{\text{char}} = 1$), then it is charging at least for the minimum number of required time slots $V_e$. The parameter $V_e$ represents the number of time slots, starting from and including time slot of event $e$, needed to meet the minimum charging time requirement.

$$\sum_{j\in J}\sum_{n\in N_j}\sum_{e^*=e}^{e+V_e-1}\delta_{j,k,n,e^*}^{dis} \geq \tau^m\lambda_{k,e}^{dis}, \forall k \in K, e \in E \quad (27)$$

Constraint (27) makes sure that if bus $k$ starts discharging at event $e$ (i.e., $\lambda_{k,e}^{dis} = 1$), then its discharging duration is at least the minimum discharging time $\tau^m$.

$$\sum_{j\in J}\sum_{n\in N_j}\sum_{e^*=e}^{e+V_e-1} y_{j,k,n,e^*} \geq V_e\lambda_{k,e}^{dis}, \forall k \in K, e \in E \quad (28)$$

Constraint (28) ensures that if bus $k$ starts discharging at event $e$ (i.e., $\lambda_{k,e}^{dis} = 1$), then it is discharging at least for the minimum number of required time slots $V_e$. The parameter $V_e$ represents the number of time slots, starting from and including time slot of event $e$, needed to meet the minimum discharging time requirement.

$$\sum_{e=\epsilon_k^g}^{\epsilon^f}\lambda_{k,e}^{char} \leq 1, \forall k \in K \quad (29)$$

Constraint (29) ensures that each bus can only connect for overnight charging once. This constraint was introduced to reflect a more realistic operational scenario. In practice, PTOs often permit buses to connect for overnight charging only once, due to practical considerations and operational constraints such as the need for manual intervention (Abdelwahed et al., 2024).

**Constraints related to PV power generation**

$$\mu_{j,e} \leq Q_{j,e}, \forall j \in J, e \in E \quad (30)$$

$$\mu_{j,e} \leq \sum_{n\in N_j}\sum_{k\in K}\frac{\delta_{j,k,n,e}^{char}}{60}\eta_{j,n}^{char}\alpha_{j,n}, \forall j \in J, e \in E \quad (31)$$

Constraint (30) and constraint (31) enforce that the usage of solar PV electricity at bus depot *j* during time slot of event $e$ when the PV panels are capable of generating electricity does not surpass the minimum of the solar energy supply $Q_{j,e}$ at bus depot $j$ during time slot of event $e$, and the total charging need for EBs charging at depot $j$ during time slot of event $e$.

$$\mu_{j,e} + \pi_{j,e} \leq Q_{j,e}, \forall j \in J, e \in E \tag{32}$$

Constraint (32) limits the solar PV energy storing $\pi_{j,e}$ and usage $\mu_{j,e}$ within the solar PV yield $Q_{j,e}$ occurring at bus depot $j$ during time slot of event $e$.

$$z_{j,e} \leq (1 - \theta_e) \sum_{n \in N_j} \sum_{k \in K} \frac{\delta^{char}_{j,k,n,e}}{60} \eta^{char}_{j,n} \alpha_{j,n}, \forall j \in J, e \in E \tag{33}$$

Constraint (33) prevents the electricity in ESS of being used for charging EBs when the PV panels are capable of generating electricity (i.e., $\theta_e = 1$). It also ensures that the electricity taken from the ESS when the PV panels do not generate electricity (i.e., $\theta_e = 0$) should not exceed the EBs' charging demand.

Linear constraints to calculate $\lambda^{\text{char}}$ and $\lambda^{\text{dis}}$ are described by (34)-(38).

$$v^{\text{char}}_{k,e} \leq \sum_{j \in J} \sum_{n \in N_j} x_{j,k,n,e-1}, \forall k \in K, e \in \{2,3, \dots, |E|\} \tag{34}$$

$$v^{\text{char}}_{k,e} \leq \sum_{j \in J} \sum_{n \in N_j} x_{j,k,n,e}, \forall k \in K, e \in E \tag{35}$$

$$v^{\text{char}}_{k,e} \geq \sum_{j \in J} \sum_{n \in N_j} \left( x_{j,k,n,e} + x_{j,k,n,e-1} \right) - 1, \forall k \in K, e \in \{2,3, \dots, |E|\} \tag{36}$$

$$\lambda^{\text{char}}_{k,e} = \sum_{j \in J} \sum_{n \in N_j} x_{j,k,n,e} - v^{\text{char}}_{k,e}, \forall k \in K, e \in E \tag{37}$$

$$v^{\text{char}}_{k,1} = 0, \forall k \in K \tag{38}$$

Constraints (34)-(38) define the variable $\lambda^{char}_{k,e}$, which indicates whether a bus $k$ starts charging at the occurrence of event $e$. To achieve this, the model uses the decision variable $v^{\text{char}}_{k,e}$ that equals 1 only if bus $k$ is charging during time slot of event $e - 1$. If this is the case, then $\lambda^{char}_{k,e}$ will be equal to 0 , meaning charging is ongoing and not starting. If $v^{char}_{k,e}$ equals 0 , then $\lambda^{char}_{k,e}$ will be equal to 1 only if bus $k$ starts charging at the occurrence of event $e$ (i.e., was not charging during time slot of event $e - 1$).

$$v^{dis}_{k,e} \leq \sum_{j \in J} \sum_{n \in N_j} y_{j,k,n,e-1}, \forall k \in K, e \in \{2,3, \dots, |E|\} \tag{39}$$

$$v^{dis}_{k,e} \leq \sum_{j \in J} \sum_{n \in N_j} y_{j,k,n,e}, \forall k \in K, e \in E \tag{40}$$

$$v^{dis}_{k,e} \geq \sum_{j \in J} \sum_{n \in N_j} \left( y_{j,k,n,e} + y_{j,k,n,e-1} \right) - 1, \forall k \in K, e \in \{2,3, \dots, |E|\} \tag{41}$$

$$\lambda^{dis}_{k,e} = \sum_{j \in J} \sum_{n \in N_j} y_{j,k,n,e} - v^{dis}_{k,e}, \forall k \in K, e \in E \tag{42}$$

$$v^{dis}_{k,1} = 0, \forall k \in K \tag{43}$$

Constraints (39)-(43) define the variable $\lambda^{dis}_{k,e}$ which indicates whether a bus $k$ starts discharging at the occurrence event $e$. To achieve this, the model uses the decision variable $v^{\text{dis}}_{k,e}$ that equals 1 only if bus $k$ is discharging during time slot of event $e - 1$. If this is the case, then $\lambda^{dis}_{k,e}$ will be equal to 0 , meaning discharging is ongoing and not starting. If $v^{dis}_{k,e}$ equals 0 , then $\lambda^{\text{dis}}_{k,e}$ will be equal to 1 only if bus $k$ starts discharging at the occurrence of event $e$ (i.e., was not discharging during time slot of event $e - 1$).

Linear constraints to calculate $a^{\mathrm{char}}$ and $a^{\mathrm{dis}}$ are as follows.

$$a_{k,e}^{char} = \sum_{j \in J} \sum_{n \in N_j} x_{j,k,n,e} - v_{k,e+1}^{char}, \forall k \in K, e \in E \setminus \epsilon^f \tag{44}$$

$$a_{k,e}^{char} + x_{j,k,n,e+1} \geq x_{j,k,n,e}, \forall k \in K, e \in E \setminus \epsilon^f, j \in J, n \in N_j \tag{45}$$

$$a_{k,e}^{dis} = \sum_{j \in J} \sum_{n \in N_j} y_{j,k,n,e} - v_{k,e+1}^{,dis}, \forall k \in K, e \in E \setminus \epsilon^f \tag{46}$$

$$a_{k,e}^{dis} + y_{j,k,n,e+1} \geq y_{j,k,n,e}, \forall k \in K, e \in E \setminus \epsilon^f, j \in J, n \in N_j \tag{47}$$

$$a_{k,\epsilon^f}^{char} = 1, \forall k \in K \tag{48}$$

$$a_{k,\epsilon^f}^{dis} = 1, \forall k \in K \tag{49}$$

Constraints (44)-(49) define the binary variables $a^{\mathrm{char}}$ and $a^{\mathrm{dis}}$ which represent if bus $k$ stops charging or discharging, respectively, during time slot of event $e$. These variables also ensure that once a bus begins charging or discharging on a specific charger, it continues using the same charger for the entire charging or discharging period. The linkage is established by using the decision variables $v_{k,e}^{char}$ and $v_{k,e}^{dis}$ which connect the charger assignment variables $x_{j,k,n,e}$ and $y_{j,k,n,e}$ to the corresponding stop indicators $a_{k,e}^{char}$ and $a_{k,e}^{dis}$.

**Restriction on arrival or departure of bus and start of charging moment**

$$\lambda_{k,e}^{dis} + \lambda_{k,e}^{char} \leq l_{j,k,e} M_{j,e}, \forall k \in K, e \in E, j \in J \tag{50}$$

Constraint (50) restricts bus $k$ to only start charging or discharging at station $j$ at the occurrence of event $e$ if there is an arrival or departure at that station (i.e., $M_{j,e} = 1$).

**Peak Power Constraints**

$$\sum_{l \in L} u_l^{\mathrm{peak}} = 1 \tag{51}$$

Constraint (51) ensures that only exactly one power level is taken as peak power level.

$$\sum_{j \in J} \sum_{n \in N_j} \sum_{k \in K} \eta_{j,n}^{char} \alpha_{j,n} x_{j,k,n,e} - \sum_{j \in J} 60 \frac{\mu_{j,e}}{\Delta T_e} - \sum_{j \in J} 60 \frac{z_{j,e}}{\Delta T_e} \leq \sum_{l \in L} U_l^{price} u_l^{\mathrm{peak}}, \forall e \in E \tag{52}$$

Constraint (52) defines the peak power level. This approach makes the assumption that the energy drawn from the ESS or from the PV panels is constant across the duration of a time slot. This assumption may or may not be valid depending on the duration of the time slot. In the worst case, however, it will be overly conservative because $\Delta T_e$ will always be larger or equal than the actual duration of the period where PV and ESS energy will be used during time slot of event $e$. And so, the peak power (and thus the final costs) might be an overestimation but never an underestimation of the actual peak power costs.

$$\sum_{j \in J} \sum_{n \in N_j} \sum_{k \in K} \left( \eta_{j,n}^{char} \alpha_{j,n} x_{j,k,n,e} - \frac{1}{\eta_{j,n}^{dis}} \beta_{j,n} y_{j,k,n,e} \right) \leq U^{\max}, \forall e \in E \tag{53}$$

Constraint (53) imposes an upper limit on the total power used to mitigate grid issues and maintain a reliable charging environment.

**Definition of Decision Variables**

$$x_{j,k,n,e}, \lambda^{char}_{k,e}, v^{char}_{k,e}, y_{j,k,n,e}, \lambda^{dis}_{k,e}, v^{dis}_{k,e}, u^{peak}_{l}, a^{char}_{k,e}, a^{dis}_{k,e} \in \{0,1\}, \forall j \in J, n \in N_j, k \in K, e \in E \tag{54}$$

$$\delta^{char}_{k,e}, \delta^{dis}_{k,e}, e_{k,e}, w^{+}_{j,e}, w^{-}_{j,e}, d_{k,e}, \mu_{j,e}, z_{j,e}, \pi_{j,e}, H_{j,e}, m^{anc}_{j,e} \geq 0, \forall j \in J, n \in N_j, k \in K, e \in E \tag{55}$$

Constraints (54) and (55) define the binary and continuous variables, respectively.

## 4. Case study

We applied the proposed DEO-based model to a case study designed to closely reflect real-world operational conditions. To enhance the practical relevance of the analysis, a potential collaboration with STIB-MIVB, the public transport operator in Brussels, was explored with the objective of facilitating data sharing and strengthening the practical relevance of the analysis. While available data was limited, sufficient high-quality open data exists for Line 53, a route that is already partially operated using electric buses. This setting therefore constitutes a realistic and well-documented test case for validating the proposed integrated charging optimization framework under real operational constraints.

The choice to focus on a single bus line is intentional and serves as a controlled validation setting. It allows the effects of integrated charging/discharging, on-site PV generation, ESS, and V2G operation decisions to be isolated and analyzed in detail, without confounding interactions across multiple routes. This choice aims to establish model validity and interpretability at the line level before extending to larger networks. The proposed formulation is not tailored to a specific bus line and can be directly extended to multiple lines or full-city bus networks without structural modifications. Incorporating multiple bus lines or an entire city bus network is the next step in this work. The scalability of the proposed approach is demonstrated via a set of test instances in Section 5.3. All data used in the model correspond to the year 2023, the most recent year for which complete datasets are available. A Julia extension in a Jupyter Notebook was used to develop the code with support from OpenAI's ChatGPT. The Gurobi 12.0.1 solver was used to formulate and solve the models, with a computational time limit of 3 hours or a relative MIP gap of 1 percent as stopping conditions. The experiments were conducted on a desktop equipped with an Intel(R) Core(TM) i7 processor with 4 threads and 16 GB of RAM. The data is freely available on https://github.com/tym2021/Integrated-charging-scheduling-for-electric-buses-with-V2G-and-renewables.git

### 4.1. Study area and bus fleet operational characteristics

The schedule of line 53 in Brussels was obtained from the STIB-MIVB website (https://transitapp.com/en/region/brussels/stib/bus-53). Line 53 runs from Westland Shopping to Hopital Militaire and vice versa. The STIB-MIVB recently placed an order for 36 eCitaros articulated electric buses; thus, the case study uses this model's characteristics to define battery capacity, lifecycle, and replacement cost[1]. The minimum and maximum allowed SOC were set as 20% and 85%, respectively (Xu et al., 2018). The initial SOC and the SOC at the end of the operations' day were set with practical considerations in mind, and both were set to 50% for testing the model. The initial SOC, however, is a parameter that would vary in real-world operations, as it depends on the battery level at the end of the previous day. In practice, this value would be set daily by the PTO. The fleet size and the number of trips for line 53 are 28 and 232, respectively. These parameters, along with other EB and ESS characteristics, are detailed in Table 2.

Table 2. ESS and EB characteristics.

| ESS Characteristics | | EB Characteristics | |
|---|---|---|---|
| $\overline{H}_j$ | 1228 kWh | Fleet size | 28 |

[1] https://www.mercedes-benz-bus.com/fr_LU/models/ecitaro/technology/battery-technology.htm

| $SOC^{min}$ | 20% | Number of trips | 232 |
|---|---|---|---|
| | | $C_k^{bat}$ | 491 kWh |
| | | $E_k^{min}$ | 25% |
| | | $E_k^{max}$ | 85% |
| | | $E_k^0, E_k^{end}$ | 50% |
| | | $N_k^{cyc}$ | 4000 |
| | | $R_k$ | 128.47€/kWh |
| | | $\tau^m$ | 5 minutes |

The closest overnight depot to this line is the Marly depot (close to the terminal Hopital Militaire), for which estimates of the additional time and distance covered during deadhead trips were calculated. The assignment of buses to the different trips was also available on the schedule. We made slight adjustments to render the assignment viable. The energy consumption of buses for each trip was estimated using the method described by Fiori et al. (2021). Note that more elaborate energy consumption models can be used by considering various influential factors, such as traffic, weather, and passenger load (Ma et al., 2021).

$$e = 0.01005v^2 - 0.3113v + 3.484 \tag{56}$$

where $e$ is energy consumed (in kWh), and $v$ is the average speed on the trip (m/s). The average speed is calculated for each trip because, despite following the same route, conditions vary significantly throughout the day due to factors such as traffic congestion. The average speed and energy consumption across all trips are 4.63 m/s (i.e., 16.68 km/h) and 17.71 kWh, respectively. Distances were available from the 2023 STIB-MIVB statistics report (STIB-MIVB, 2023), with slight adjustments for the deadhead trips that were added when it was a trip that arrived from or returned to the Marly depot.

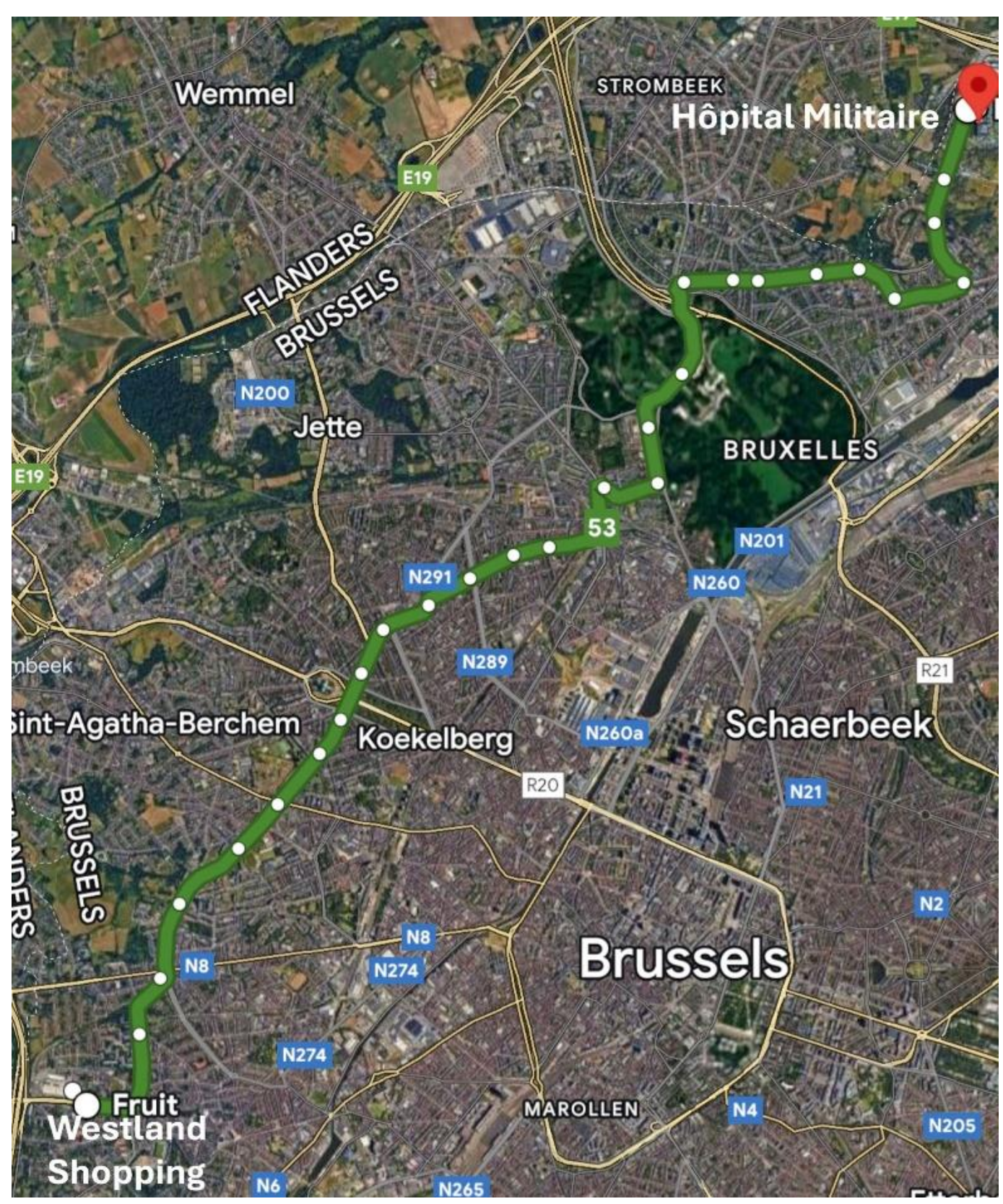

Figure 3. Service area of Line 53.

4.2. Charging infrastructure at the overnight depot and terminals

We consider two types of charging infrastructure: depot charging and opportunity charging. The depot charging is used when buses return to the depot at the end of a service day for overnight charging. Charging at the depot typically uses plug-in technology operating at relatively low power levels, ranging from 30 to 150 kW (Ashkezari et al., 2024). Fast chargers are located at the line terminals, allowing buses to recharge during layover between two trips. The charging technology in this setting is a pantograph, which usually provides power levels between 150 kW and 600 kW (Ashkezari et al., 2024). Values of charger powers were chosen as 150 kW for overnight chargers and 300 kW for terminal chargers, consistent with industry standards. For their efficiency, a relatively high performance level was assumed, as STIB-MIVB has been actively investing in electric infrastructure in recent years. Based on this, chargers are probably part of a modern installation, with a deployment no earlier than 2018. Efficiencies were set at 92%, a typical value for recently installed equipment. This parameter is expected to decrease over the charger's lifetime. The parameters used for the charging infrastructure specification is shown in Table 3. With regard to the choice of the number of chargers at each charging station, we tune this parameter in order to find a balance between the model's feasibility and complexity.

Table 3. Charger specifications by location.

| Location | Charger | Charging efficiency | Charging power (kW) | Discharging efficiency | Discharging power (kW) |
|---|---|---|---|---|---|
| Marly Depot (j=1) | 1 | 0.92 | 150 | 0.92 | 120 |
| | 2 | 0.92 | 150 | 0.92 | 120 |
| | 3 | 0.92 | 150 | 0.92 | 120 |
| Terminal Westland Shopping (j=2) | 1 | 0.92 | 300 | 0.92 | 240 |
| | 2 | 0.92 | 300 | 0.92 | 240 |
| Terminal Hopital Militaire (j=3) | 1 | 0.92 | 300 | 0.92 | 240 |
| | 2 | 0.92 | 300 | 0.92 | 240 |

4.3. Renewable energy integration

To obtain information on the solar radiation in Brussels, an open dataset made available by the Royal Meteorological Institute of Belgium was used. It contains hourly data points for solar radiation in 2023, and, as with electricity prices, the average hourly solar radiation values were computed to obtain aggregate values. The values can be found in Table A.1 in the appendix. In practice, the PTO would provide day-ahead market prices and solar radiation forecasts as inputs to the model for scheduling the next day's charging operations. Solar energy generation was calculated based on an assumed PV surface area of 1876.6 $m^2$ at the overnight depot. This figure comes from the STIB-MIVB 2023 Statistics Report, which states that 9383 $m^2$ of solar panels are distributed across five depots (STIB-MIVB, 2023). Assuming an equal distribution, the solar yield during a charging event was calculated by multiplying the solar radiation per square meter by the PV area, then adjusting for event duration and converting the result to kWh. We assume that there are no PV panels present at the fast-charging terminal stations because the surface available to install PV panels is limited, and bus terminals are often situated in urban locations where buildings and trees create shadowing during the day. For the battery energy storage system, we assume that only the Marly Depot is equipped with one with a capacity of 1228 kWh.

4.4. Peak power and electricity pricing

STIB-MIVB manages its tram and metro systems in collaboration with the Brussels distribution system operator, Sibelga. At the Marly depot, the electric bus charging infrastructure is connected to the existing metro power network, which operates at medium voltage. Based on this setup, the case study assumes similar grid connection conditions and related capacity charges, estimated at approximately 49.5 €/kW per year. In this framework, the peak power cost is calculated by multiplying the highest daily peak power drawn by the daily tariff. These values are based on the 2023 electricity tariff document published by Brugel (Brugel, 2022). For charging terminals, medium-voltage connections are assumed. This is justified by the fact that a single 300 kW charger already exceeds the safe operational limits of typical low-voltage infrastructure. When two or three chargers are installed at the same terminal, the total site load typically necessitates a medium-voltage connection (11 kV- 25 kV), especially when reliability and power quality are critical operational concerns (O'Regan, 2024). To model peak power charges, a trial-and-error approach was adopted to define ten discrete operating bands ranging from 100 kW to 1000 kW, with $U^{max}$ set at 1000 kW. Table A.2 in the appendix lists the corresponding cost for each peak power level.

Regarding the electricity price, it plays a central role in cost optimization. Belgian DA electricity prices were used for 2023 from the ENTSO-E platform to build a "typical day" profile by averaging hourly prices across the year. The values obtained can be found in Table A.1 in the appendix. To calculate electricity selling prices, the case study assumes a conservative approach of 25% lower than the purchasing price, which aligns with market practices (Manzolli et al., 2025).

### 4.5. Scenario description

This section walks through the sequence of scenarios (i.e., the DEO model variants with different components) evaluated. It starts with a bare-bones baseline—no peak-power charges, no V2G capability, and no on-site renewables or ESS. Then, we design a scenario where peak power costs, V2G, and degradation costs are included. Finally, a scenario with all extensions (i.e., peak power, V2G, degradation cost, PV panels, and ESS) is designed to look into the interactions between all components.

#### 4.5.1. Basic Scenario

This scenario optimizes charging operations based on day-ahead prices of electricity, with the objective of minimizing daily operational costs associated with charging events. All other features of the model are disabled: peak power costs, V2G, battery degradation costs, and renewable energy production with energy storage solution are excluded from both the objective function and the constraints. The discussion of results subsequently computes peak power costs incurred in this scenario by examining the number of chargers used simultaneously during each event.

#### 4.5.2. Scenario considering peak power costs and V2G with degradation costs

The second scenario includes V2G and degradation costs. This is the only modification compared to the previous scenario. The goal is to assess the profitability of selling energy back into the grid while accounting for the additional battery capacity fade caused by discharging activities. To enhance the computational efficiency, discharging is restricted to occur only during predefined peak hours which are set from 7:00 to 10:00 and from 18:00 to 21:00. This assumption is based on the premise that PTOs will only sell electricity back to the grid when market prices are sufficiently high, typically during peak hours when electricity prices reach their daily maximum. These constraints narrow the solution space for $y_{j,k,n,e}$ , thereby accelerating the resolution of this model.

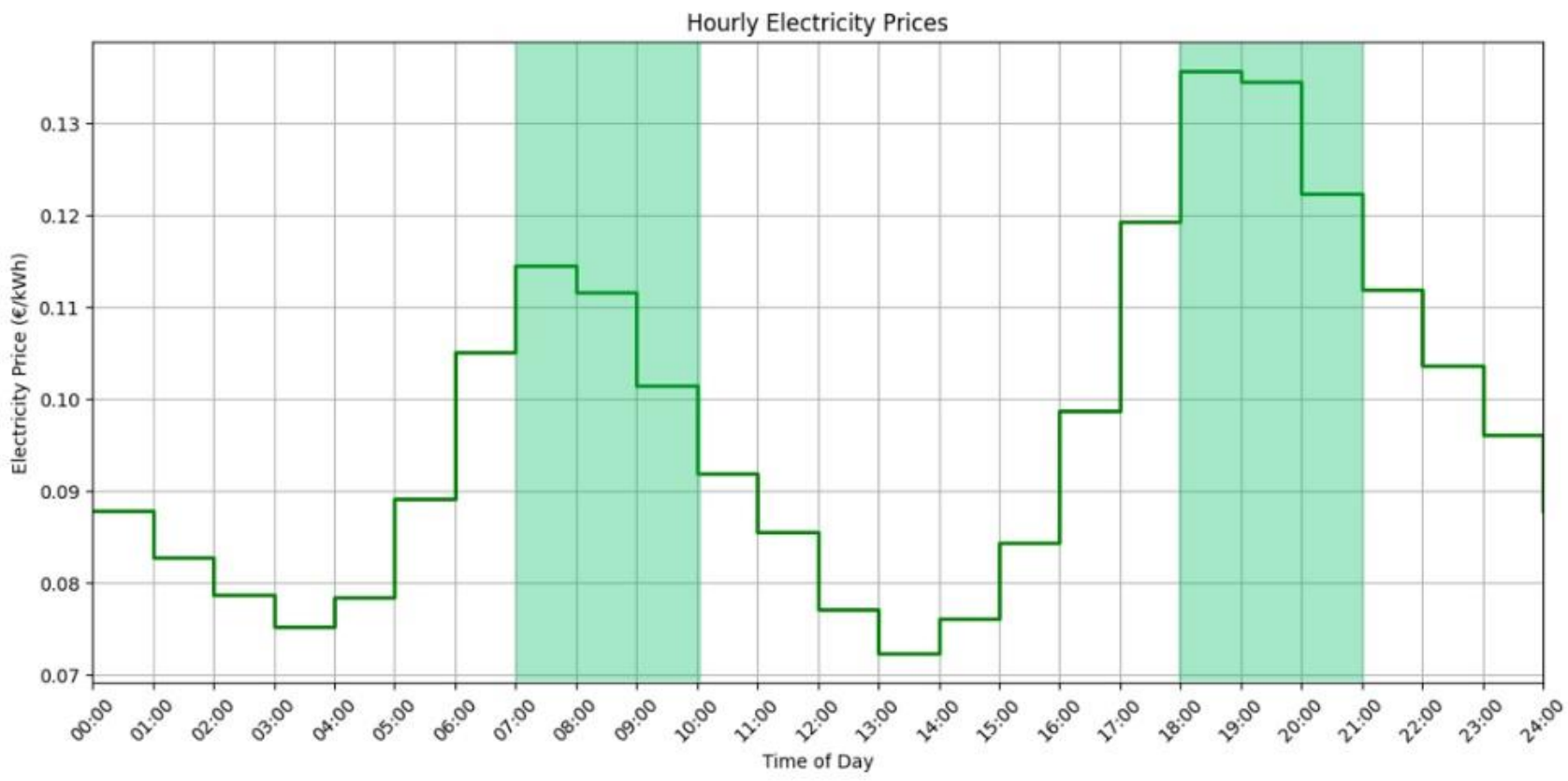

Figure 4. Evolution of hourly electricity prices over a "typical" day. Shaded areas indicate the peak hours during which discharging events are allowed.

4.5.3. Scenario with all extensions

This final scenario incorporates all the previously discussed features. In this scenario, electricity generated from local PV panels installed at the overnight depot can be used to charge either the EBs or the ESS at zero cost during periods of solar availability. Stored electricity can also be drawn from the ESS at zero cost for electric buses during events where there is no sun. Furthermore, the ESS can participate in the energy market by re-injecting electricity into the grid when necessary, providing an additional revenue stream for the PTO through energy arbitrage. The objective is to analyze the individual impact of each component on the case study, as well as its combined effects when integrated into a unified system. Table 4 provides a comprehensive comparison of the different scenarios.

Table 4. Overview of scenarios with their characteristics.

| Scenario | Peak Power | V2G | Degradation Costs | PV and ESS |
|---|---|---|---|---|
| Basic | | | | |
| PP + V2G + DC | ✓ | ✓ | ✓ | |
| All | ✓ | ✓ | ✓ | ✓ |

Remark: PP: Peak Power. DC: Degradation Costs

## 5. Results

This section provides a comprehensive analysis of the computational results obtained. Comparison of the scenarios is done in three steps: first, an analysis of the cost implications of the different extensions; next, an examination of how they reshape the charging schedule; and finally, an assessment of each scenario's effect on energy production and consumption. Based on these results, managerial insights and decision implications are derived.

5.1. Costs analysis

Table 5 represents the solutions obtained for each scenario with detailed decomposition of costs, including charging costs, discharging revenues, battery degradation costs, and peak power costs, as well as the maximum peak power reached, the MIP gap, and the CPU time. The result shows that our full model can be solved with 1% optimality gap within 2 hours for the problem size of 28 V2G-enabled EBs and 232 trips, demonstrating the computational efficiency of our DEO-based MILP formulation. The total operating costs for the full model (scenario all) is significantly lower (359.9 euros) compared to the benchmark (more than 800 euros). We analyse and compare the total costs and their cost components for these scenarios in the following subsections.

Table 5. Overview of results obtained for each scenario.

| Scenario | Total Costs (€) | Charging Costs (€) | Discharging Revenues (€) | Degradation Costs (€) | Peak Power Costs (€) | Peak Power (kW) | MIP Gap | CPU Time (s) |
|---|---|---|---|---|---|---|---|---|
| Basic | 859.58 | 691.66 | 0.00 | 0.00 | 167.92 | 1242 | 1.00% | 2443 |
| PP+V2G+DC | 813.61 | 718.97 | 0.00 | 0.00 | 94.64 | 700 | 1.00% | 11679 |
| All | 359.90 | 476.35 | 197.58 | 0.00 | 81.12 | 600 | 0.99% | 6931 |

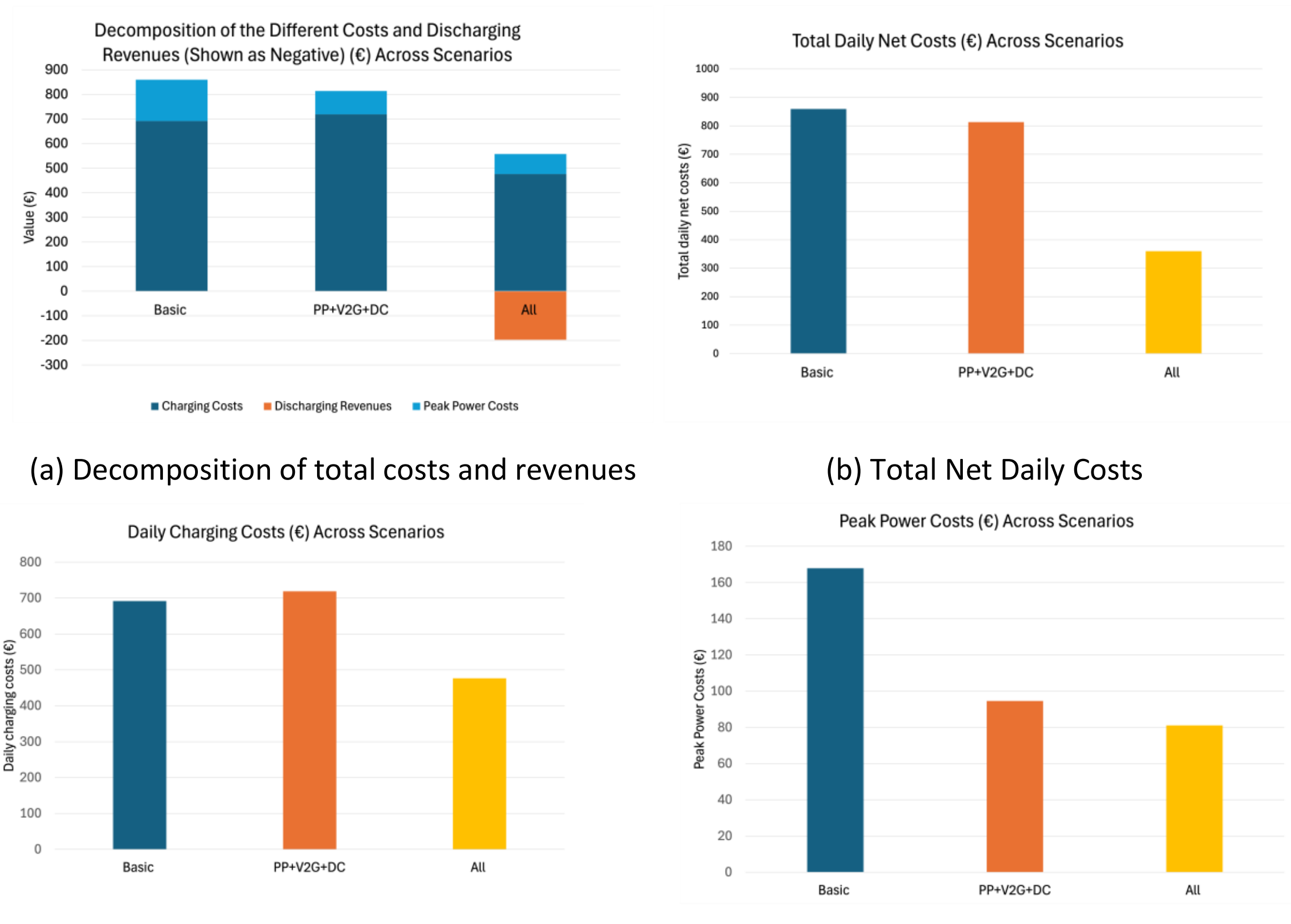

(a) Decomposition of total costs and revenues

(b) Total Net Daily Costs

(c) Daily Charging Costs

(d) Daily Peak Power Costs

Figure 5. Comparison of results between the different scenarios.

*Total net cost analysis*

From the decomposition of costs in Figure 5a, it is observed that including peak power constraints decreases the total cost by 5% compared to the baseline scenario. The cost reduction between the Basic and PP+V2G+DC scenarios is solely due to the inclusion of peak power charges in the objective function, as V2G is not utilized, its degradation costs outweigh the benefits of discharging. The final cost settles at €813.60 in this scenario, reflecting an attempt to optimize net charging costs and peak power costs while preserving battery health. Then, the introduction of solar panels and a battery energy storage system, starting from the PP+V2G+DC scenario, drastically reduces the final net cost, from €813.61 to €359.90, a reduction of approximately 56% (see Figure 5b).

**Managerial implication**: The optimal cost obtained in the All scenario represents a decrease of roughly 58% compared to the basic scenario, highlighting the importance of using an optimization tool and extensions to reduce daily operating costs for the PTO. These results raise interesting implications considering the future evolutions of the potential of V2G. In particular, it suggests that its inclusion may not be economically meaningful today, but could become so in the future, depending on the evolution of battery replacement costs, its degradation behaviour, and tariff margin.

*Charging costs analysis*
Figure 5c illustrates the daily charging costs. The results reveal that the inclusion of peak power in the model increases charging costs by 4%. This is due to the model's dual objective: minimizing charging expenses while also accounting for peak demand, a factor that, when all costs components are considered, emerges as a significant contributor to overall operational costs. As a result, buses are sometimes charged during periods with higher electricity prices to avoid triggering high peak power charges. The optimized (dis)charging operations balances energy price minimization against peak demand exposure. The introduction of solar panels and an ESS has a significant impact, resulting in a drastic reduction of 34% in total charging costs due to local PV meeting the charging demand at zero cost.

**Managerial implication**: The result confirms that peak power acts as a binding operational constraint, which reshapes charging schedules rather than merely imposing an additional cost. The ESS further amplifies this effect by coupling PV generation and enabling surplus solar energy to be stored or used during high electricity price periods. This demonstrates that the integrated PV–ESS configuration with demand charges can significantly lower the charging costs and reduce exposure to energy price volatility and peak demand charges.

*Peak power cost analysis*
The introduction of peak power constraints and associated costs in the model leads to a noticeable reduction in maximum peak power drawn from 1242 kW in the basic scenario to 700 kW when peak power pricing is considered. This strategic limitation of simultaneous charging leads to a peak power cost reduction of approximately 44% (see Figure 5d). In the baseline scenario, peak power costs represent up to 20% of total net costs. Including peak power constraints reduces this share to 11%, a notable 9% drop. Further improvements are observed with the integration of solar panels and an energy storage solution: the maximum peak demand decreases from 700 kW to 600 kW, resulting in an additional peak power cost reduction of approximately 14%. This is consistent with the expected benefit of relieving grid stress, as part of the demand is offset by local renewable generation and stored energy. Interestingly, scenarios that include V2G and battery degradation costs show similar levels of peak power and associated costs as the Peak Power scenario, suggesting that these features alone do not significantly influence peak reduction when compared to the impact of solar and storage integration.

**Managerial implication**: The result highlights the importance of effective peak power management for cost-effective charging management. In particular, it implies that PTOs can achieve significant charging cost savings without changes to charging infrastructure by coordinating fleet charging scheduling. The integration of PV generation and ESS can enable further savings for peak power demand charges and enhance grid stability.

*Discharging revenue and degradation costs analysis*

The introduction of solar panels, and more significantly, the integration of the ESS, results in a substantial increase in discharging revenues, rising from zero to €197.58, as seen in Table 5. This increase is primarily driven by the ESS performing energy arbitrage: storing energy when prices are low and selling it back to the grid during high-price periods. Figure 6 illustrates this mechanism, showing how ESS exports align with high electricity prices. This highlights the model's ability to exploit price fluctuations for economic gain. In the Basic scenario and the PP+V2G+DC scenario, there are no discharging revenues. Further tests will be conducted to explore these aspects in more depth. It is also worth noting that a limitation of the current model is the absence of degradation costs for the ESS, which could lead to a slightly overestimation of its economic contribution. Regarding the degradation costs, they are never observed in the scenarios studied. This is because, by introducing degradation costs into the objective function, the revenue that can be generated with V2G is not large enough to offset the additional battery wear incurred. The sensitivity analysis of the battery degradation costs will be presented in the next section.

**Managerial implication**: The results show that integrating V2G capability alone is insufficient for generating charging/discharging arbitrage gains without an ESS. This implies that V2G-enabled charging and discharging flexibility must be coupled with an ESS to generate higher revenue during periods of high electricity prices. The results also underscore that V2G profitability depends on its battery degradation costs.

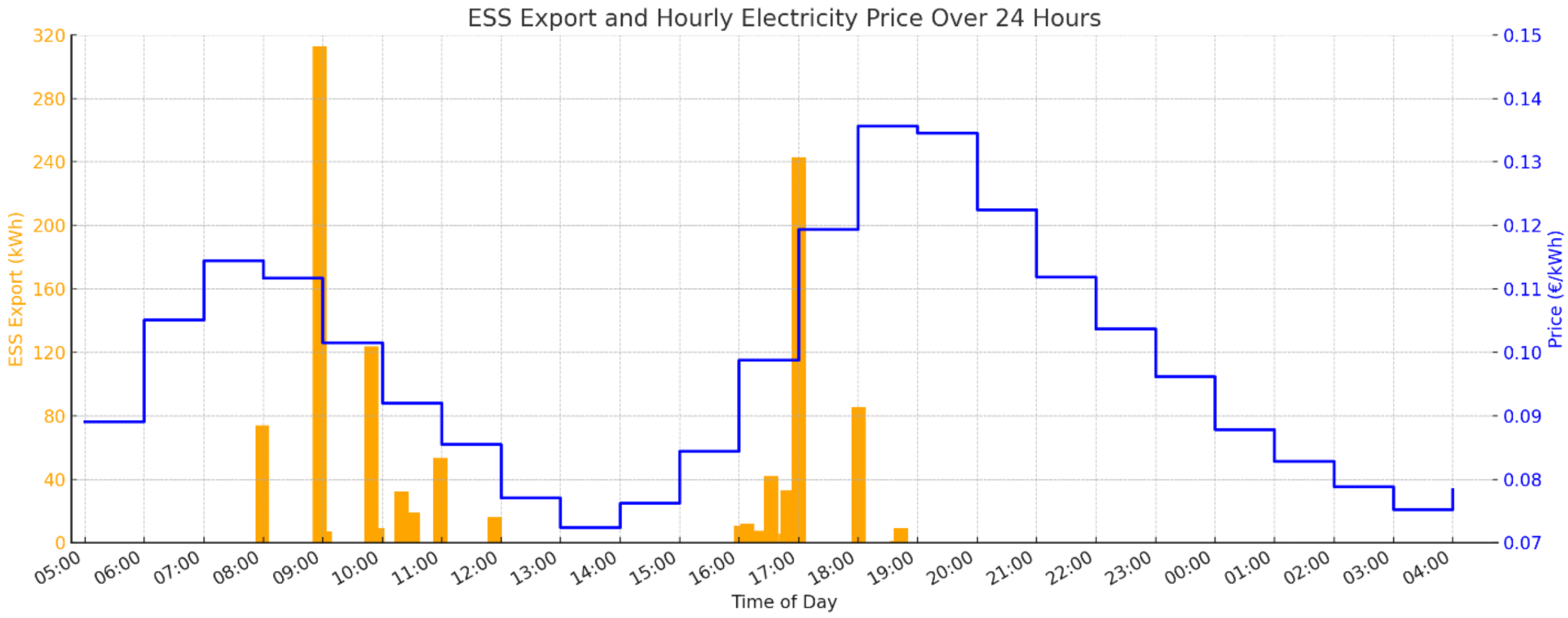

Figure 6. ESS Export and hourly electricity price over 24 hours.

*Conclusion*

Taken together, the scenario results point to a clear hierarchy of cost drivers. Peak-power charges are the first lever: just setting a limit on how much power the depot can pull from the grid at once cuts the total electricity bill by about 5%. Then, V2G, on today's tariff margin and replacement-cost assumptions and with internalized battery wear, is not used. By contrast, adding on-site photovoltaics and an energy-storage system reshapes the cost picture. The PV+ESS system wipes out roughly one-third of the baseline expenditure (i.e., charging cost and demand charges), and permits an additional revenue stream through energy arbitrage. Of the €359.90 per-day cost that remains, about 77% is the net result of energy bought minus discharging revenue, while the remaining 23% comes from peak-power charges. This hypothesis will be verified in the sensitivity analysis.

These findings suggest a clear prioritization strategy for PTOs: Prioritize PV and ESS deployment and peak shaving under current conditions for sustained cost reduction. V2G might become more appealing if battery prices decrease or export margins increase significantly, depending on future market trends.

## 5.2. Analysis of energy produced & consumed

Different operational strategies influence energy sourcing and exportation. In the Basic scenario and the PP+V2G+DC scenario, all energy (8113 kWh) is imported from the grid due to the absence of alternative energy sources. The energy landscape changes significantly with the addition of a photovoltaic array and a stationary battery. Grid dependency drops markedly, from 8113 kWh in the initial setup to 5623 kWh in the All scenario, which represents a substantial reduction of approximately 30%. Further analysis shows that, on this representative day selected for the case study, the solar panels produce up to 4818 kWh. From that output, 2409 kWh (≈50 %) directly charge the buses, 80 kWh (≈ 1.5%) pass through the ESS before returning to the fleet, about 170 kWh remain stored within the ESS (≈ 3.5%), and 2159 kWh (≈ 45%) are cycled through the ESS and exported to the grid for additional revenue. The latter refers to energy arbitrage, where surplus solar energy is fed back into the grid to generate additional revenue for the PTO. In other words, almost the entire ESS throughput, about 96%, is devoted to revenue-earning export rather than to vehicle charging. This effect is discernible in the diagram presented in Figure 7, where the energy flows of the All scenario are represented. Overall, the combination of renewable energy, V2G, and energy storage not only reduces grid imports but also transforms the PTO into an active contributor to the broader energy network.

**Managerial implication**: A key managerial insight is that nearly all ESS throughput is used for exporting energy to the grid rather than directly charging buses. This suggests that the primary economic value of ESS lies in arbitrage opportunities rather than in supporting vehicle charging alone. The result illustrates how V2G and ESS can be leveraged without negatively affecting bus availability while generating additional revenue.

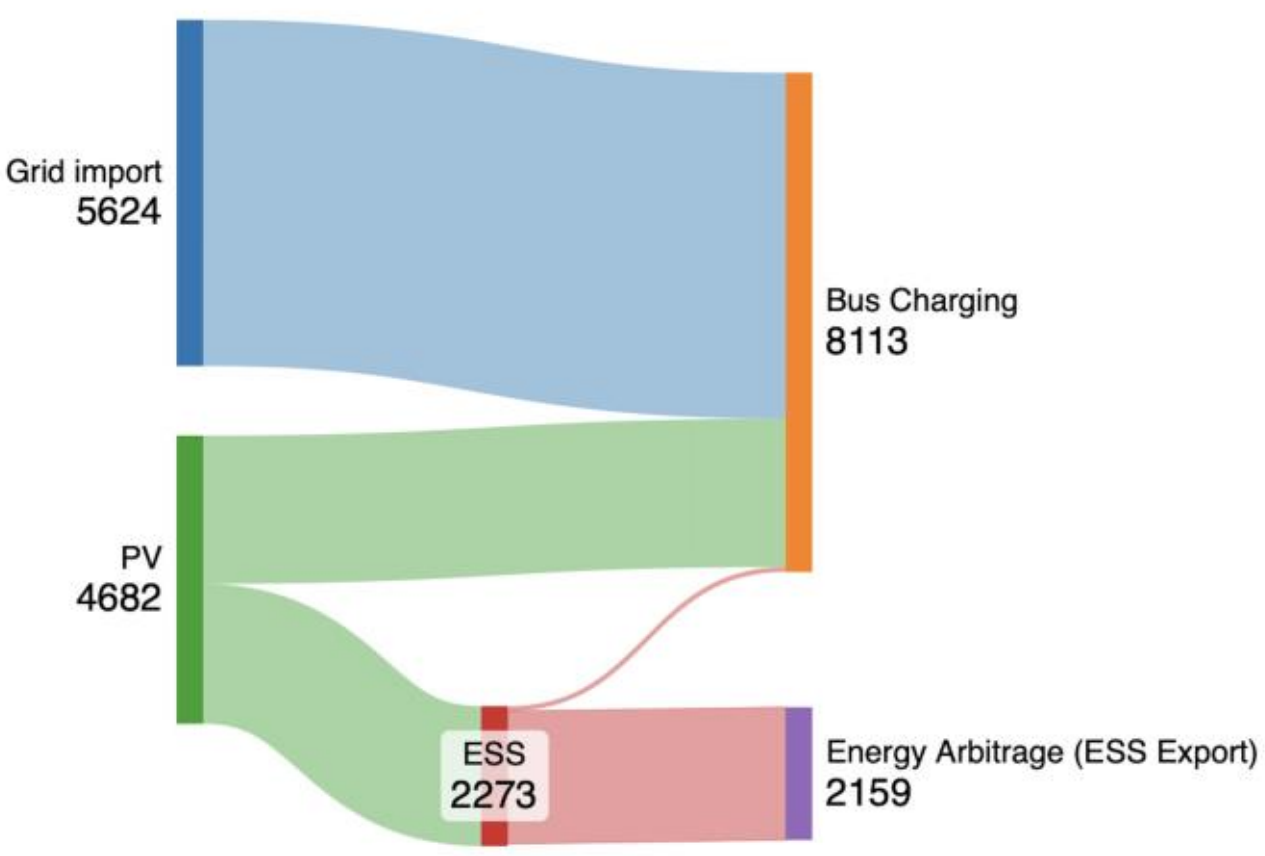

Figure 7. Energy flows (kWh) of the all scenario.

5.3. Scalability and computational performance

In this section, we investigate how changes in the number of buses, the capacity of the ESS storage, the efficiency of chargers, and the PV panels' surface affect the computational results. We generate 16 test instances by varying the above parameters. In the results shown next, instance names are labeled as follows: “B” number of buses; “E” charging efficiency; “PV” PV panels surface; “ESS” ESS storage capacity. For example, the instance labelled as 10B-82E-938.3PV-614ESS corresponds to a situation where there are 10 buses running on the line that can charge at chargers that have 82% efficiency with the surface of PV panels available equal to 938.3 m$^2$ and the ESS having a storage capacity of 614 kWh.

The results are shown in Table 6. This result comparison across instances with varying parameter values unveils some interesting observations. Firstly, computational times are shorter for scenarios with 10 buses. This is not surprising, as fewer buses reduce the model’s complexity by decreasing the number of decision variables and constraints. The reduced charging demand also leads to lower peak power and charging costs. Most notably, the instance with 10 buses, high charging efficiency, and the maximum capacity for both PV panels and ESS results in a negative net total cost. This indicates a tipping point at which the PTO generates net revenues, emphasizing the financial viability of high renewable integration under favorable conditions. Furthermore, reducing PV by half often results in a more than 50% drop in discharging revenues, which significantly worsens net charging costs and suggests nonlinear benefits from solar investment.

Table 6. Cost breakdown and performance metrics by configuration.

| Configuration | Total Costs (€) | Charging Costs (€) | Discharging Revenues (€) | Degradation Costs (€) | Peak Power Costs (€) | MIP Gap (%) | CPU Time (s) |
|---|---|---|---|---|---|---|---|
| 28B-82E-1876PV-614ESS | 379.25 | 483.24 | 180.06 | 8.46 | 67.60 | 2.62 | 10799 |
| 28B-82E-938PV-1228ESS | 570.99 | 554.70 | 51.30 | 0.00 | 67.60 | 0.56 | 10801 |
| 28B-82E-938PV-614ESS | 580.06 | 528.30 | 29.36 | 0.00 | 81.12 | 2.00 | 10800 |
| 28B-92E-1838PV-1228ESS | 357.93 | 483.45 | 193.12 | 0.00 | 67.60 | 1.34 | 10801 |
| 28B-92E-1876PV-614ESS | 380.61 | 479.48 | 180.00 | 0.00 | 81.12 | 1.86 | 15469 |
| 28B-92E-938PV-1228ESS | 569.43 | 563.60 | 75.29 | 0.00 | 81.12 | 0.94 | 6184 |
| 28B-92E-938PV-614ESS | 575.85 | 553.89 | 59.16 | 0.00 | 81.12 | 1.66 | 10801 |
| 28B-92E-1876PV-1228ESS | 359.90 | 476.36 | 197.58 | 0.00 | 81.12 | 0.99 | 6931 |
| 10B-82E-1876PV-614ESS | 25.52 | 258.48 | 273.52 | 0.00 | 40.56 | 1.00 | 10514 |
| 10B-82E-938PV-1228ESS | 175.58 | 282.85 | 147.83 | 0.00 | 40.56 | 1.00 | 921 |
| 10B-82E-938PV-614ESS | 194.07 | 275.32 | 121.81 | 0.00 | 40.56 | 0.99 | 10077 |
| 10B-92E-1876PV-1228ESS | 2.39 | 258.48 | 296.64 | 0.00 | 40.56 | 0.65 | 3597 |
| 10B-92E-1876PV-614ESS | 21.33 | 251.91 | 271.14 | 0.00 | 40.56 | 1.27 | 2273 |
| 10B-92E-938PV-1228ESS | 175.78 | 285.76 | 150.54 | 0.00 | 40.56 | 1.00 | 3615 |
| 10B-92E-938PV-614ESS | 194.42 | 285.85 | 132.00 | 0.00 | 40.56 | 0.97 | 1987 |
| 10B-92E-1876PV-1228ESS | -1.79 | 251.91 | 294.26 | 0.00 | 40.56 | 0.00 | 2650 |

Remark: The CPU time is limited to 3 hours, or until a relative MIP gap of 1% is reached, whichever occurs first.

## 6. Sensitivity analysis

As discussed, V2G's anticipated commercial readiness within the next few years justifies its inclusion in this study. Consequently, the scenarios assume that the STIB-MIVB can export electricity to the grid both from ESS coming from renewable energies and through V2G from electric buses. In the results' analysis, the model currently places limitations on the use of V2G. This may be due to overly restrictive parameters. Therefore, the sensitivity analysis examines the potential impact of reduced battery costs and explores how changes in the coefficient of the wholesale price of electricity affect outcomes. This includes not only the tariff margin but also the influence of a previously unconsidered factor: green certificates. The sensitivity analyses are all conducted on the All scenario (full model) with 28 buses using the parameter setting described in the sections 4.1-4.4.

### 6.1. Battery degradation cost

The battery cost values are taken from a study by Mauler et al. (2021). This study aggregates findings from a variety of different sources published between 2010 and 2020, and it offers a comprehensive picture of forecasted values from 2035 to 2050. The sensitivity analysis was performed for the years and values outlined in Table 7. Since the study provided values for the costs in $, a projected equilibrium exchange rate of 1.10 $/€ is applied, reflecting a scenario in which global financial conditions remain subdued due to ongoing geopolitical tensions and economic uncertainty (Martínez et al., 2025). The analysis shows that the use of V2G increases when the battery cost decreases. The V2G integration starts to appear in the results in 2040, when the battery costs are equal to 101.20 €/kWh. The use of V2G, however, is minimal in absolute size. In 2050, the amount generated when selling back to the grid through V2G is €6.32. This number only represents 1.12% of the total costs (i.e., charging costs, peak power costs, and degradation costs), which is a tiny portion. For more detailed information about the optimization results, actual values are available in Table 8.

**Managerial implication**: Despite the small values for discharging revenues, technological innovations in the field of lithium-ion batteries have a positive impact on the amount of V2G undertaken. A reduction in battery replacement costs could provide an attractive new revenue stream for PTOs in the future if combined with other factors such as advantageous electricity selling prices or an increase in the total energy that the batteries of buses can provide throughout their lifespan.

Table 7. Projected battery costs by year (Mauler et al., 2021).

| Year | Projected Battery Cost (€/kWh) |
| --- | --- |
| 2023 | 128.47 |
| 2035 | 119.90 |
| 2040 | 101.20 |
| 2045 | 88.00 |
| 2050 | 78.10 |

Table 8. Sensitivity analysis of battery cost.

| Year | Total Net Costs (€) | Degradation Costs (€) | V2G Revenue (€) | MIP Gap (%) | CPU Times (s) |
| --- | --- | --- | --- | --- | --- |
| 2023 | 359.90 | 0.00 | 0.00 | 0.99 | 6644 |
| 2035 | 360.00 | 0.00 | 0.00 | 1.00 | 5520 |
| 2040 | 360.96 | 1.04 | 5.51 | 1.00 | 3418 |
| 2045 | 360.33 | 1.13 | 5.81 | 1.00 | 4829 |
| 2050 | 360.92 | 1.10 | 6.32 | 0.99 | 5570 |

### 6.2. Tariff margin and green certificates

Based on the work of Manzolli et al. (2025), small producers of renewable energy currently receive about 75% of the wholesale market price for energy they inject into the grid. This discount is primarily due to the need for small producers to sell through aggregators or energy suppliers, as direct market access requires a pooled capacity of at least 1 MW. These intermediaries introduce transaction costs that reduce the effective price received. In Brussels, energy producers also benefit from **green certificates (GCs)**, which provide supplementary income to support renewable energy generation. This mechanism sets Brussels apart from Flanders and Wallonia, where similar support schemes have already been phased out (Bobex, n.d.). For example, small PV installations (< 5 kW) receive approximately 2.05 green certificates per MWh. In the STIB-MIVB case, with an evaluated installed capacity between 281.49 kW and 319.022 kW, the system falls into the large installation category (i.e., >250 kW), entitling it to 0.58 GC per MWh in 2025 (Brugel, 2025). The full tariff table is available in Table A.3 in the appendix and on the official website. Each green certificate was traded between €65 and €90 in 2023 in Belgium, hovering around the same values today. This implies that STIB-MIVB could receive an additional €38 to €52 per MWh on top of the energy sale price. When combined with the base energy injection revenue (i.e., around 75% of the wholesale price), the total return from grid injection can exceed 100% of the wholesale market rate, making export economically viable under current conditions. Brussels maintains this certificate mechanism to meet its renewable electricity target of at least 12% by 2025 and to continue incentivizing solar energy production (International Energy Agency, 2018). However, as installation costs fall and renewable generation scales up, the region is expected to either phase out this support or recalibrate it, reflecting developments already seen in other Belgian regions. Green certificates were excluded from the initial analysis under the assumption that the program had been retired. However, after confirming its continued application, it is now included in the sensitivity analysis to account for any short-term financial impact it may still have for public transport operators like STIB-MIVB. While not critical to long-term planning due to its likely phase-out, the mechanism currently provides a valuable revenue stream that shapes near-term energy dispatch and export decisions.

### 6.2.1. Scenario description for the wholesale price evolution

We consider three scenarios related to the projection of the percentage of the wholesale price that bus operators could receive for injected electricity from 2025 to 2050:

- **Positive scenario (structural enablers)**: We assume that policies increasingly empower prosumers and aggregators, enabling easier and more direct market participation for smaller players (EU Commission, 2024). This improved access allows operators like STIB-MIVB to capture a greater share of the wholesale market value without relying heavily on intermediaries, from 75% in 2025 to 100% in 2050, for the energy price recovered from exports. The overall policy and market environment remain favorable, encouraging investment in V2G and solar infrastructure.
- **Conservative scenario (market saturation and price pressure)**: The conservative scenario assumes modest policy support and growing solar penetration without adequate market reform. As solar renewable capacity increases, midday oversupply becomes common, reducing energy prices during peak production hours. In this context, STIB-MIVB faces diminishing returns from injected energy. While V2G and ESS allow some flexibility in shifting energy use to higher-value hours, they cannot fully counteract the structural dip in selling prices. Green certificates follow the same declining pattern as in the positive scenario, but offer less strategic value given the market constraints
- **Pessimistic scenario (policy retrenchment and structural barriers)**: In the pessimistic scenario, regulatory and market barriers are even worse, leading to a progressive decline in the fraction of the wholesale price received by bus operators. This could occur due to oversupply, the imposition of new grid charges on distributed energy injections, or specific curtailment rules that prevent producers from selling at preferred times (Elia). These combined forces further depress injection margins, stabilizing around 40% by 2050. In this negative scenario, by 2050, bus operators are on their own in the market and still face structural disadvantages, capturing significantly less than the wholesale price for their energy.

For each scenario, the energy price portion, which is the percentage of the wholesale price achieved through energy sales, and the green certificate support, which is a bonus expressed as a percentage of the wholesale price, are shown in Table A.4.2 and Table A.4.3 in the appendix, respectively. The total combined value is summarized in Table A.4.1 in the appendix.

### 6.2.2. Results

We solve the DEO model for three scenarios. The obtained results are analysed as follows. Figure 8 shows the different energy flow changes with respect to different tariff margins. It reveals a distinct change in operating behavior around the 75-80% tariff-capture threshold. Below this point, the system prioritizes self-consumption, a higher proportion of the photovoltaic output is routed straight into the buses, with the grid covering any shortfall. The ESS supports this strategy because, for example, when the tariff-capture rate is at 60% of the wholesale price, approximately 45% of its discharged energy is allocated to bus charging. However, once the margin reaches about 80%, exporting becomes more profitable. The control strategy, therefore, flips: buses are increasingly recharged from the grid, and the energy produced by the PV is steered towards export in order to capture the higher selling price. As a result, the share of PV energy used directly for vehicle charging decreases, and ESS discharging is repurposed for generation rather than bus charging. For instance, at an 80% tariff margin, only 12 kWh of ESS output is allocated to bus charging, and from 90% onward, this share falls to 0%. Notably, the first instance of V2G export (≈ 15 kWh) appears around 80%, disappears briefly, and then re-emerges around 95% with 45 kWh, reflecting a point where export revenues begin to sufficiently offset degradation costs. Then, V2G exports increase sharply from 45 kWh to 1766 kWh as the margin rises from 95% to 110%. This surge is accompanied by a continued increase in ESS exports and grid imports, as the system exploits low off-peak prices to recharge and sell back at higher tariffs —a classic arbitrage response. This leads to a rise in degradation costs from just €1.58 at 95% to over €61 at 110%. This shows the trade-off between short-term revenue gains and long-term battery wear. As shown in Figure 9a, at a 110% margin, charging costs are offset mainly by discharging revenues, with remaining costs primarily driven by battery degradation and peak power demand charges. Furthermore, Figure 9b confirms the significant economic impact of changes in the tariff margin: daily net charging costs decrease by approximately 42%, dropping from €414 at a 40% margin to €238 at 110%. The complete set of results for the different tariff margins can be found in Table A.5 in the appendix.

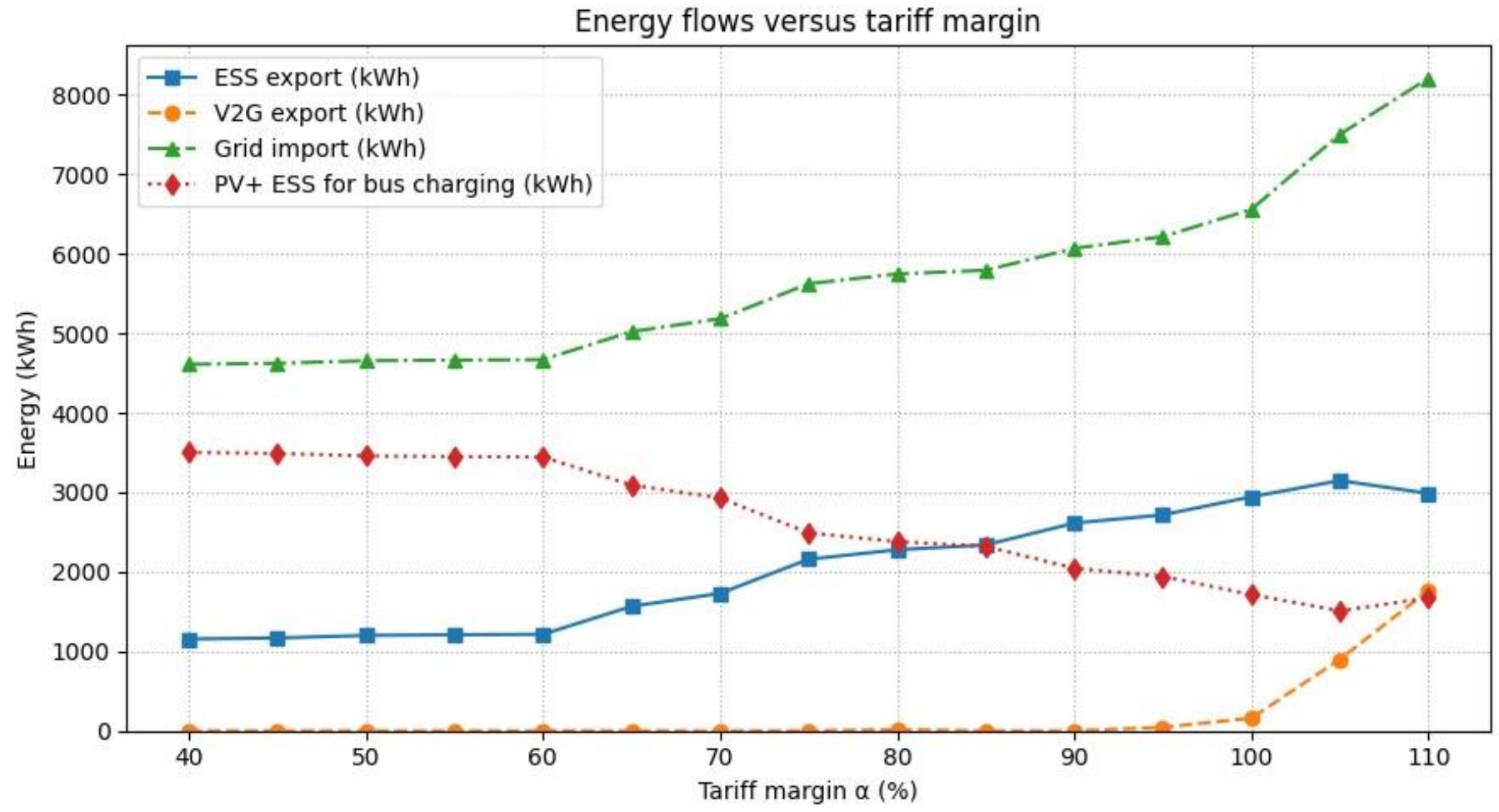

Figure 8. Energy flows versus tariff margin.

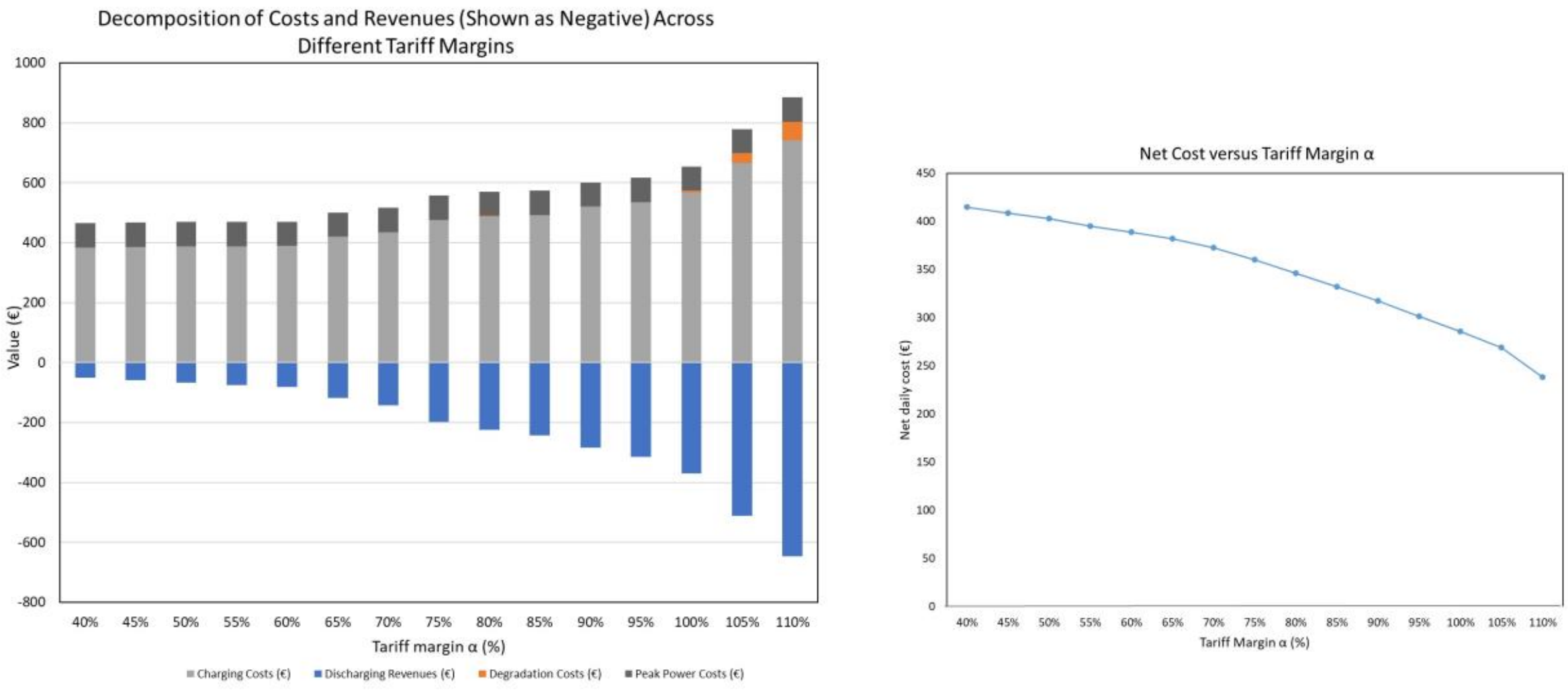

(a) Decomposition of costs and revenues (b) Daily charging cost versus tariff margin

Figure 9. Sensitivity analysis results on charging costs and cost decomposition.

By considering the different scenarios and incorporating green certificate incentives, clear trends are observed in the total net daily cost over time (see Figure 10). In all scenarios, current conditions favor bidirectional charging, as the present 110% tariff margin favors V2G through the financial support provided by green certificates. Figure 11 shows that, currently, GCs cover up to 35% of potential discharging revenues, offering a significant advantage to PTOs. However, this advantage is temporary. As the scheme is progressively phased out across all scenarios, the analysis uncovers a general upward trend in net daily costs, particularly in the conservative and pessimistic cases.

In the conservative and pessimistic scenarios, while the model attempts to shift as much solar production as possible to high-value periods using intelligent charging strategies, the declining injection margins and oversupply pressure prevent significant gains. These scenarios remain heavily dependent on ESS for flexibility, with V2G relatively inactive due to weak market signals. Essentially, they become "ESS-only" strategies.

Only in the positive scenario does the combination of favorable market access, tariff margin growth allow meaningful V2G monetization. By 2050, this scenario achieves 30% lower net operating costs than the others, primarily by capitalizing on discharging revenues through ESS and V2G once tariff margins exceed 95% and degradation penalties become manageable.

**Managerial implication**: Given these findings, the PTO is strongly incentivized to act now: both to capitalize on the short-term gains from green certificate support and to position itself for medium-term benefits by investing in additional ESS capacity. Even in less favorable scenarios, a well-sized ESS offers resilience against volatile grid prices and future market restrictions. In favorable conditions, it also serves as a revenue-generating asset through arbitrage opportunities.

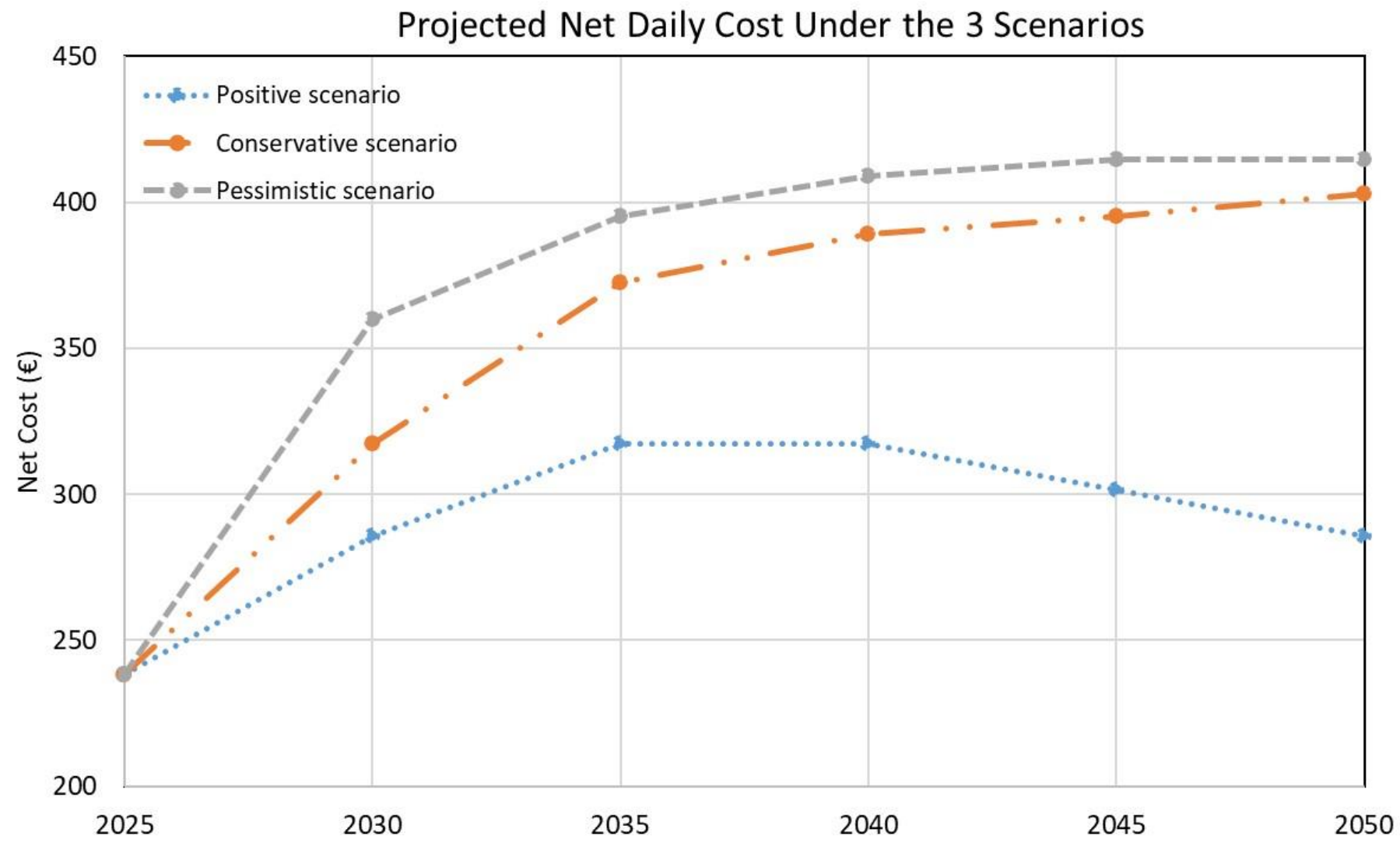

Figure 10. Projected total net daily costs under the different scenarios.

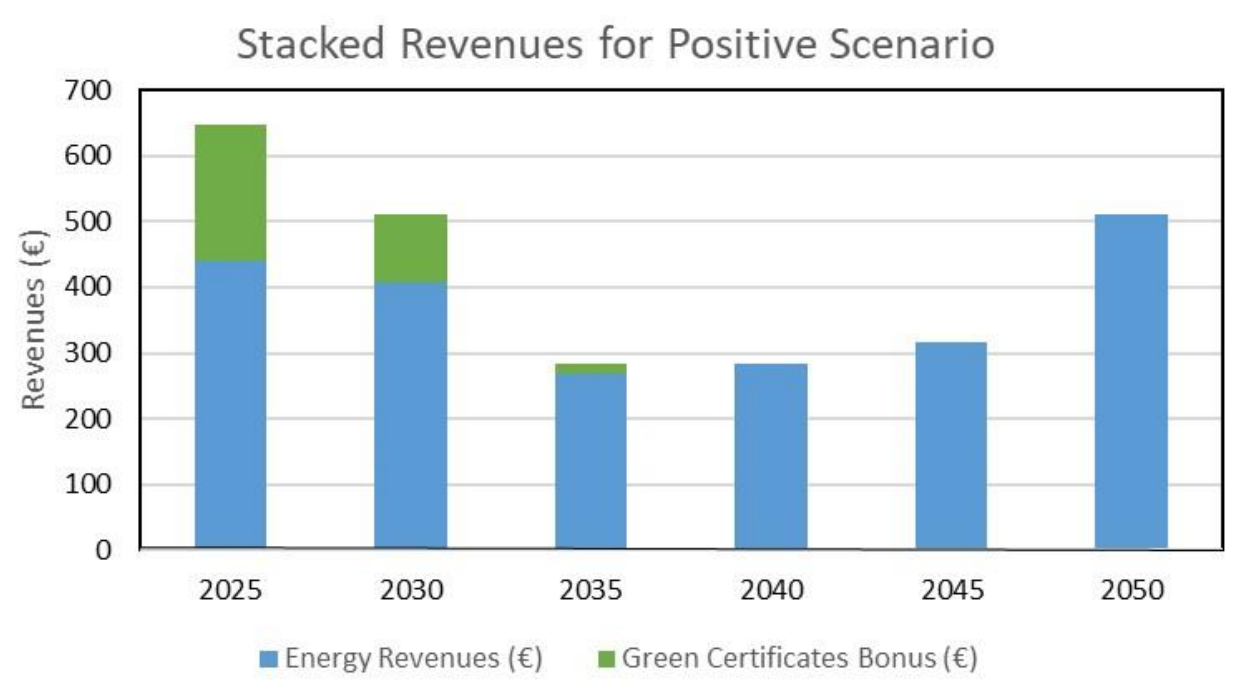

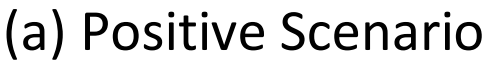
(a) Positive Scenario

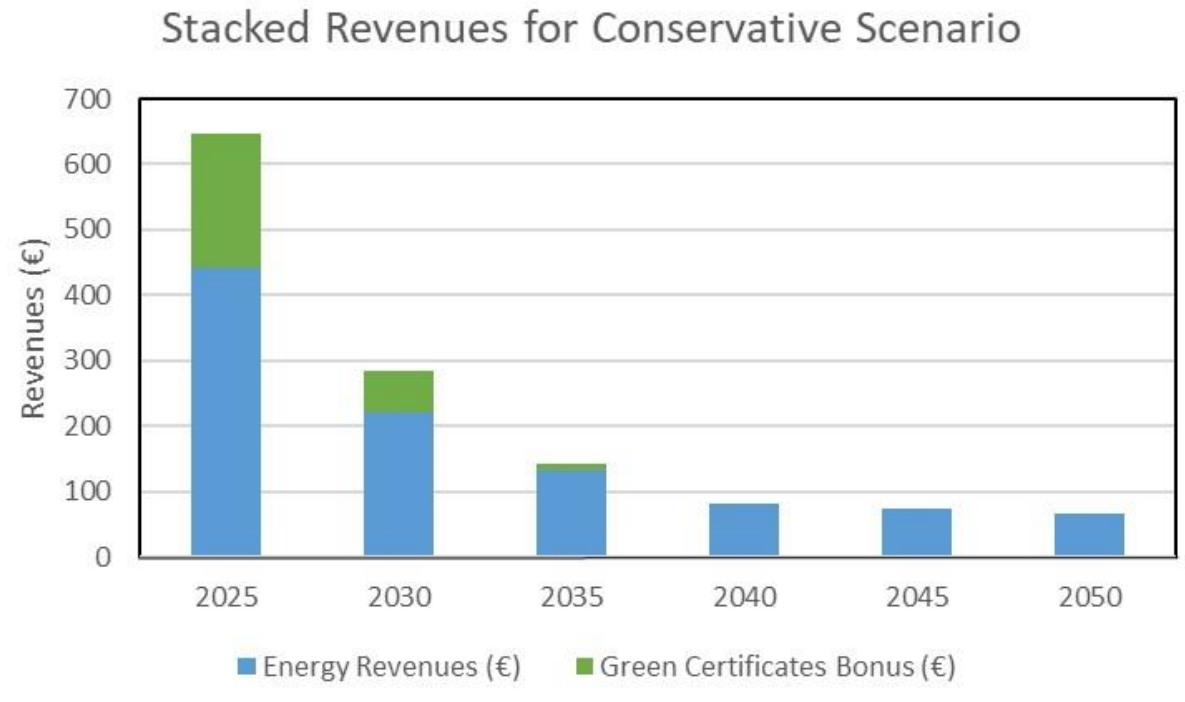

(b) Conservative Scenario

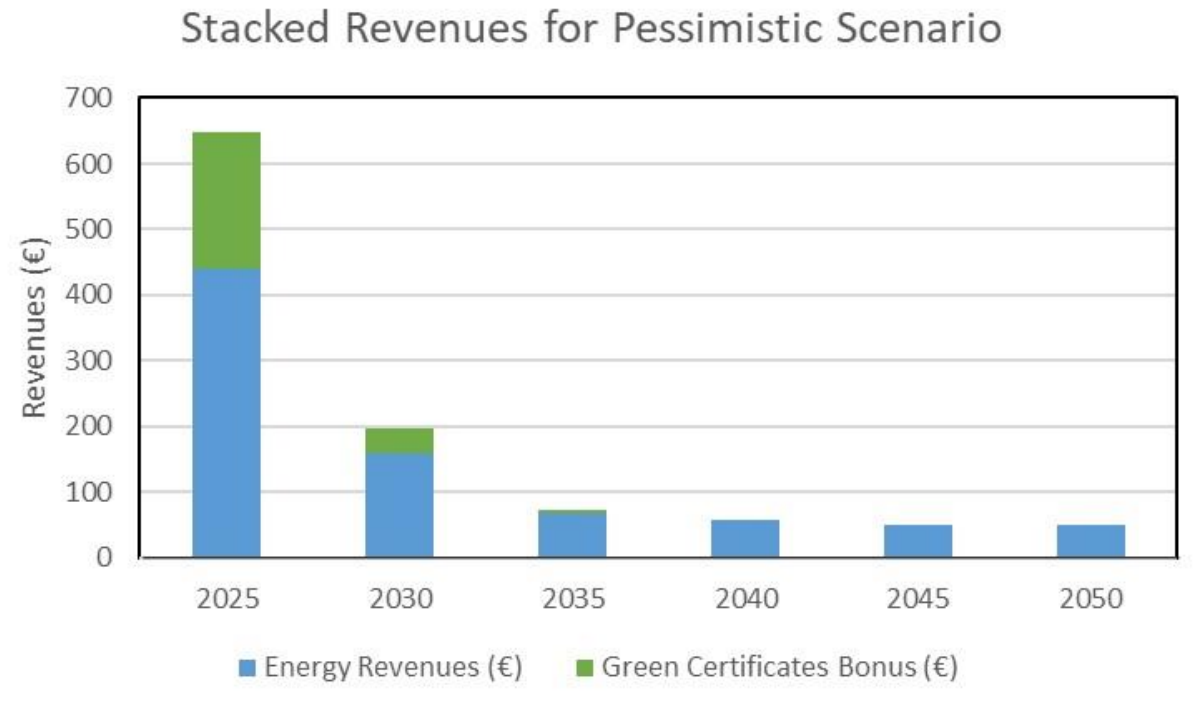

(c) Pessimistic Scenario

Figure 11. Stacked revenues for each scenario.

**7. Conclusions and discussion**

This study develops a comprehensive modeling framework for coordinating the daily charging and discharging operations of an electric bus fleet. The primary methodological contribution lies in the adoption of an event-driven time discretization approach to represent complex charging and discharging decisions as well as the interactions among multiple system components and the power grid. In contrast to conventional discrete-time discretization schemes, the proposed event-driven approach reflects practical operational constraints and substantially enhances computational efficiency and scalability for real-world applications, while preserving the temporal resolution required for accurate modeling of operational decisions.

While existing studies typically incorporate dynamic energy prices, capacity-based demand charges, on-site PV generation, ESS, V2G capabilities, and battery degradation costs in isolation or in limited combinations, the proposed model integrates all of them within a scalable optimization framework that enables cost-effective EB charging decisions and supports grid stability. The case study results demonstrate that such integration is essential for achieving flexible and economically efficient charging management, as it allows the model to capture and resolve the structural trade-offs between charging, energy storage, and grid interaction decisions.

We apply the model for a real-world case study involving 28 electric buses and 232 trips operated by STIB-MIVB in Brussels on a single route over a 24-hour period. The results show that the proposed DEO model provides the lowest net operational cost with a total decrease of 58% compared to the benchmark without these extended components. The result demonstrates the significant value of using an optimization tool that considers all extensions for PTOs. We also test the scalability and computational performance of the model on a set of test instances. The results demonstrate the computational efficiency of our approach, as near-optimal solutions could be obtained by the state-of-the-art commercial MIP solver within 2 hours, with a 1% optimality gap. In addition, we assessed the contribution of each extension to the system's performance and conducted sensitivity analyses to examine how degradation cost and tariff margin influence the extent to which V2G is utilized. Several managerial insights are drawn to support operational and planning decisions of PTOs.

These findings imply that public transport operators must not only plan for infrastructure upgrades but also consider smart energy strategies such as leveraging dynamic electricity pricing, demand charges, integrating local renewable energy and ESS, and deploying flexible technologies like V2G to remain both economically viable and environmentally sustainable. While the model delivers useful insights and demonstrates strong potential for practical application, a few limitations remain that could influence its performance under real-world conditions and should be considered for future refinement. Firstly, although the model includes battery degradation costs for electric buses, it does not incorporate degradation effects on the ESS. Secondly, the electricity pricing framework relies on day-ahead prices and does not capture additional market opportunities such as the intraday market, which may affect pricing accuracy. Finally, external uncertainties affecting the energy consumption of buses (e.g., weather variability, driver behavior) were not included. Ignoring these factors may reduce the robustness of the proposed scheduling strategy. This study contributes to the growing body of research helping PTOs design cost-effective, flexible, and sustainable charging strategies that capitalize on renewable integration, storage systems, and smart grid opportunities such as V2G and peak power management.

Future research could enhance the model by incorporating stochastic elements to account for energy consumption uncertainty and by expanding the number of actors in the landscape to reflect the complex electricity pricing structure that exists beyond the day-ahead market scheme. Since detailed data on multiple routes and their associated buses were not available, the model could not be scaled to a full network of routes and depots. However, doing so would be an important step in assessing how the model performs in a more operationally realistic context. Finally, explicitly modeling ESS degradation would provide valuable operational insights for public transport operators and support more informed energy storage strategies.

**Appendix**

Table A.1. Average energy purchasing price and solar radiation by hour (2023)

| Time | Price (€/kWh) | Solar Radiation (W/m$^2$) |
|---:|---:|---:|
| 00:00 – 1:00 | 0.0877 | 0.00 |
| 1:00 – 2:00 | 0.0827 | 0.00 |
| 2:00 – 3:00 | 0.0776 | 0.00 |
| 3:00 – 4:00 | 0.0752 | 0.00 |
| 4:00 – 5:00 | 0.0780 | 0.26 |
| 5:00 – 6:00 | 0.0890 | 10.20 |
| 6:00 – 7:00 | 0.1050 | 53.27 |
| 7:00 – 8:00 | 0.1144 | 113.66 |
| 8:00 – 9:00 | 0.1116 | 173.29 |
| 9:00 – 10:00 | 0.1015 | 226.68 |
| 10:00 – 11:00 | 0.0919 | 267.26 |
| 11:00 – 12:00 | 0.0855 | 255.47 |
| 12:00 – 13:00 | 0.0771 | 241.97 |
| 13:00 – 14:00 | 0.0724 | 249.12 |
| 14:00 – 15:00 | 0.0762 | 247.96 |
| 15:00 – 16:00 | 0.0843 | 227.26 |
| 16:00 – 17:00 | 0.0897 | 177.09 |
| 17:00 – 18:00 | 0.1193 | 129.41 |
| 18:00 – 19:00 | 0.1356 | 81.24 |
| 19:00 – 20:00 | 0.1345 | 27.62 |
| 20:00 – 21:00 | 0.1224 | 0.83 |
| 21:00 – 22:00 | 0.1118 | 0.00 |
| 22:00 – 23:00 | 0.1036 | 0.00 |
| 23:00 – 00:00 | 0.0961 | 0.00 |

Table A.2. Peak power rating and corresponding price

| Index | Peak Power (kW) | Daily Price (€) |
|:---:|:---:|---:|
| 1 | 100 | 13.52 |
| 2 | 200 | 27.04 |
| 3 | 300 | 40.56 |
| 4 | 400 | 54.08 |
| 5 | 500 | 67.60 |
| 6 | 600 | 81.12 |
| 7 | 700 | 94.64 |
| 8 | 800 | 108.16 |
| 9 | 900 | 121.68 |
| 10 | 1000 | 135.21 |

A.3 Brugel’s Green Certificates Ratios

The table below is in French and outlines the GC allocation rate for different energy technologies as of January 1, 2025. What is most relevant for this study is the last column, titled " Taux d’Octroi (CV/MWh)", which represents the amount of green certificates granted per MWh of electricity produced. In particular, the "Photovoltaïque" section refers to PV panel installations. The table details how the GC rate varies depending on the installation type and size. For example, small rooftop PV systems receive a higher GC rate per MWh than larger installations.

Table A.3. Brugel’s Green Certificates Ratios.

| Coefficient Multiplicateur et Taux d'Octroi | | | vers: | 1-01-2025 |
|---|---|---|---|---|
| **Installation consommant du carburant** | | | | |
| Technologie | Conditions | | Coefficient Multiplicateur | **Octroi (CV)** |
| Installation de cogénération | Biomasse | / | 1 | En fonction du rendement et de l'économie de CO2 |
| | Au gaz naturel | / | 0 | 0 |
| **Installation ne consommant pas de carburant** | | | | |
| Technologie | Conditions | | Coefficient Multiplicateur | **Taux d'Octroi (CV/MWh)** |
| Taux de base | | | 1 | 1,818 |
| Photovoltaïque | Pour les installations classiques, sur base de la puissance (kWc) | [0-5] | 1,130 | 2,055 |
| | | ]5-36] | 1,074 | 1,953 |
| | | ]36-100] | 0,559 | 1,016 |
| | | ]100-250] | 0,353 | 0,642 |
| | | > 250 | 0,319 | 0,580 |
| | pour le BIPV, sur base du type* | Skylight | 1,130 | 2,055 |
| | | Garde-corps | 1,130 | 2,055 |
| | | Brise-soleil | 1,046 | 1,902 |
| | | Façade ventilée | 1,423 | 2,587 |
| Éolienne | / | / | 1 | 1,818 |
| Hydroélectrique | / | / | | |

A.4 Different revenue and green certificate bonus as a percentage of wholesale price (2025–2050)

Table A.4.1. Table Total revenue (energy + GC) as a percentage of wholesale price (2025–2050)

| Year | Positive Scenario | Conservative Scenario | Pessimistic Scenario |
|---|---|---|---|
| 2025 | 75% | 75% | 75% |
| 2030 | 80% | 70% | 60% |
| 2035 | 85% | 65% | 50% |
| 2040 | 90% | 60% | 45% |
| 2045 | 95% | 55% | 40% |
| 2050 | 100% | 50% | 40% |

Table A.4.2. Energy injection revenue as a percentage of wholesale price (2025–2050)

| Year | Positive Scenario | Conservative Scenario | Pessimistic Scenario |
|---|---|---|---|
| 2025 | +35% | +35% | +35% |
| 2030 | +20% | +20% | +15% |
| 2035 | +5% | +5% | +5% |
| 2040 | 0% | 0% | 0% |
| 2045 | 0% | 0% | 0% |
| 2050 | 0% | 0% | 0% |

Table A.4.3. Green certificate bonus as a percentage of wholesale price (2025–2050)

| Year | Positive Scenario | Conservative Scenario | Pessimistic Scenario |
|---|---|---|---|
| 2025 | 110% | 110% | 110% |
| 2030 | 100% | 90% | 75% |
| 2035 | 90% | 70% | 55% |
| 2040 | 90% | 60% | 45% |
| 2045 | 95% | 55% | 40% |
| 2050 | 100% | 50% | 40% |

Table A.5. Results of sensitivity analysis on tariff margins.

| Margin | Total Costs (€) | Charging Costs (€) | Discharging Revenues (€) | Degradation Costs (€) | Peak Power Costs (€) | MIP Gap | CPU Times (s) |
|---|---|---|---|---|---|---|---|
| 0.40 | 414.64 | 384.44 | 50.93 | 0.00 | 81.12 | 0.01 | 2782.42 |
| 0.45 | 408.77 | 385.69 | 58.04 | 0.00 | 81.12 | 0.01 | 2498.53 |
| 0.50 | 403.08 | 388.24 | 66.29 | 0.00 | 81.12 | 0.01 | 4716.49 |
| 0.55 | 395.27 | 388.27 | 74.12 | 0.00 | 81.12 | 0.01 | 6202.44 |
| 0.60 | 388.91 | 388.89 | 81.10 | 0.00 | 81.12 | 0.01 | 4447.54 |
| 0.65 | 381.81 | 419.69 | 119.00 | 0.00 | 81.12 | 0.01 | 4150.86 |
| 0.70 | 372.43 | 434.82 | 143.52 | 0.00 | 81.12 | 0.01 | 4262.94 |
| 0.75 | 359.90 | 476.36 | 197.58 | 0.00 | 81.12 | 0.01 | 6931.34 |
| 0.80 | 346.30 | 489.11 | 224.46 | 0.53 | 81.12 | 0.01 | 5113.03 |
| 0.85 | 332.00 | 493.16 | 242.29 | 0.00 | 81.12 | 0.01 | 9373.07 |
| 0.90 | 317.35 | 520.09 | 283.87 | 0.00 | 81.12 | 0.01 | 3647.00 |
| 0.95 | 301.30 | 534.35 | 315.76 | 1.59 | 81.12 | 0.02 | 10801.02 |
| 1.00 | 285.67 | 567.98 | 369.02 | 5.59 | 81.12 | 0.03 | 10624.90 |
| 1.05 | 268.71 | 667.00 | 510.70 | 31.29 | 81.12 | 0.03 | 10603.90 |
| 1.10 | 238.25 | 742.66 | 646.97 | 61.44 | 81.12 | 0.05 | 10801.53 |